\documentclass[11pt]{amsart}
\usepackage{geometry}                
\usepackage{graphicx}
\usepackage{amssymb}
\usepackage{epstopdf}
\usepackage{subfig} 
\usepackage{color}
\usepackage{caption}
\usepackage{cancel}
\usepackage{arydshln}
\usepackage[font=small,skip=0.1pt]{caption}

\newcommand {\R} {\mathbb {R}}

\newcommand {\A} {\mathbb {A}}

\newcommand {\Pj} {\mathbb {P}}

\newcommand {\no} {\noindent}

\newcommand {\pa} {\partial}
\newcommand {\rest}[1] {|_{#1}}

\theoremstyle{definition}

\newtheorem{defn}{Definition}[section]

\newtheorem{thm}{Theorem}[section]

\newtheorem{eg}[thm]{Example}
\newtheorem{rem}[thm]{Remark}

\title{The title of my page}

\title{Spelling Rules for the Monster/Semple Tower.}
\author{Alex L. Castro}
\email{alex.castro@mat.puc-rio.br}
\address{Departamento de Matem\'atica, Pont\'ificia Universidade Cat\'olica do Rio de Janeiro. Phone: +55-21-3527-1722}
\author{Wyatt C. Howard}
\email{whoward@ucsc.edu}
\thanks{W.C.H. is the corresponding author}
\address{IMPA, Instituto Nacional de Matem\'atica Pura e Aplicada, Rio de Janeiro}
\date{\today}

\begin{document}


\maketitle

\tableofcontents
\addtocontents{toc}{\protect\setcounter{tocdepth}{1}}

\begin{abstract}

We further study the incidence relations that arise from the various subtowers, known as Baby Monsters, which exist within the $\R^{3}$-Monster Tower.  This allows us to
complete the $RVT$ class spelling rules.  We also present a method of calculating the various Baby Monsters that appear within the Monster Tower.


\end{abstract}

\vspace{.2in}

\noindent {\bf MSC} Primary: $58A30$
\

\no Secondary: $58A17$, $58K50$.
\vspace{.1in}

\noindent {\bf Keywords:} Goursat Multi-Flags, Prolongation, Semple Tower, Subtowers, Monster Tower.
\vspace{.2in}


\section{Introduction}\label{sec:intro}
This paper further investigates the various critical hyperplane relations that occur within a rank-$3$ distribution associated to each level of the $\R^{3}$-Monster Tower.
Studying these hyperplane interactions yields the complete spelling rules for the $RVT$ coding system. 
In \cite{castro2} we studied how critical hyperplanes appeared within a rank $3$ distribution that is associated to each of the points within the Monster Tower.
We used these hyperplane relations in order to classify the obits within the Monster Tower.  
Classifying these orbits enabled us to classify a certain geometric distribution known as a Goursat $n$-flag.

A \textit{Goursat $n$-flag of length $k$}, also called a \textit{Goursat Multi-Flag} for $n \geq 2$, is a distribution, say $D$, of rank $(n+1)$ sitting inside of a $(n+1) + kn$ dimensional ambient manifold.
The rank of the associated flag of distributions
$$D \hspace{.2in} \subset  \hspace{.2in} D + [D,D]\hspace{.2in} \subset  \hspace{.2in} D + [D,D] + [[D,D],[D,D]] \dots ,$$
increases by $n$ at each bracketing step.

R. Montgomery and Z. Zhitomirskii showed in \cite{mont1} that the problem of classifying up to local diffeomorphism of Goursat Flags is equivalent to the problem of classifying points within an iterated 
tower of manifolds that they called the \textit{$\R ^{2}$-Monster Tower}, also known as the \textit{Semple Tower} in algebraic geometry (\cite{sempleInvest} and \cite{bercziMap}).  Then in \cite{mont2} 
they focused on classifying the various points by partitioning them into equivalence classes known
as the the $RVT$ coding system.  This coding system served as an important tool which allowed them to completely classify the points within the $\R ^{2}$-Monster Tower.  In \cite{castro} A. Castro
and R. Montgomery took the first steps toward extending this $RVT$
coding system for the $\R ^{3}$-Monster Tower.  In the $\R ^{3}$ case there are distinguished directions, labelled by the letter $L$, which come from the intersection of the vertical and
tangency hyperplanes.  While the critical hyperplane relations are completely understood in the $\R^{3}$ case over regular, vertical, tangency, and $L$ points, it is not obvious what happens over the
new critical directions that appear over the $L$ directions.  The word \textit{over} here, and throughout the rest of the paper, refers to the rank $3$ distribution which is
``over'' each point within the Monster Tower.  In \cite{castro2} we began studying the sources of these critical hyperplanes.  Over $L$ points there exist not just one, but two
types of tangency hyperplanes and three new distinguished directions, see Figure \ref{fig:three-planes}.  These tangency hyperplanes above the $L$ directions are labeled as $T_{1}$ and $T_{2}$ and 
the new $L$ directions by $L_{j}$ for $j = 1,2,3$.  One natural question to ask is the following:
\begin{quote}
{\it What are the critical hyperplane incidence relations over each of the directions $T_{i}$ and $L_{j}$ for $i = 1,2$ and $j = 1,2,3$?}
\end{quote}

This is exactly the question we are solving in this paper, where we determine the complete $RVT$ spelling rules.
We want to examine this question because it will enable us to classify points at higher levels of the Monster Tower.  In \cite{castro2} we used a
technique called the \textit{isotropy method} that allowed us to classify not only the points within the fourth level of the $\R^{3}$-Monster Tower, but to classify points at any other level of the tower so 
long as we know how to characterize the $RVT$-classes in Kumpera-Rubin coordinates.  In \cite{castro2} we only understood the critical hyperplane configurations over $R, V, T,$ and $L$ points.  
However, after level $4$ of the Monster Tower it was not entirely clear what the spelling rules are for the various $RVT$ classes.  In particular, we pointed out how the existence of the $L$ direction 
is the source of this confusion.  This means that once we complete the spelling rules, we can then apply the isotropy method to determine how many orbits lie within any level of the 
tower.  This is the main motivation for the work presented in this paper.
We will derive the $RVT$ spelling rules for the $\R ^{3}$- Monster Tower over each of the previously unknown planes and directions listed above.  As a result, this will yield the complete
$RVT$ spelling rules.  They are as follows:
\begin{enumerate}
\item Any $RVT$ code string must begin with the letter $R$.
\item $R$ :  $R$ and $V$.
\item $V$ and $T \, (=T_{1})$ : $R, V, T,$ and $L$.
\item $L \, (= L_{1})$ and $L_{j}$ for $j = 2,3$ : $R, V, T_{i}$ for $i = 1,2$, and $L_{j}$ for $j = 1,2,3$.
\item $T_{2} \, : R, \, V, \, T_{2}, \, \text{and} \, L_{3}$.
\end{enumerate}


\begin{rem}
The colon indicates the letters which can come after the given letter(s).  For example, after the letter $R$ we can add either the letter $R$ or $V$.  The letters $T$ will just be used for
$T_{1}$ and $L$ for $L_{1}$ from this point on.  Additionally, $T_{i}$ will just refer to one of the two tangency hyperplanes and $L_{j}$ as one of the three distinguished directions.
\end{rem}

Another motivation for studying these incidence relations is understanding how the $RVT$ coding system relates to the $EKR$ coding system studied by P. Mormul (\cite{mormul1}).  The initial
relationship between these two coding systems that characterize Goursat Multi-Flags was first studied in an appendix to \cite{castro2}.
One other benefit of establishing this connection is that it helps explain the relationship between the $RVT$ code and the articulated arm system.  In \cite{pelletier2} F. Pelletier and M. Slayman 
investigate the relationship between the $RVT$ codes and the $EKR$ system for points within the fourth level of the $\R ^{3}$-Monster
Tower along with the restrictions placed upon the positions of the bars within the articulated arm system \cite{respondek}.  We discuss this line of research in more detail within the conclusion section of
the paper.
\

We also want to point out that all of our work done under the name of the Monster Tower can be done in the same setting as the Semple Tower with a base of $\A ^{3}$, affine $3$-space, and with the 
terminology of Baby Monsters replaced by Subtowers.
\

In Sections $2$ and $3$ we present the main definitions needed in order understand the statements of our main results.  We also present a few illustrative examples that will help the reader understand 
the proofs of the paper.  Section $4$ consists of the proofs of our main results.  Then in Section $5$ we provide a summary of our findings and further research directions.
\vspace{.2in}

\noindent {\bf Acknowledgements.} Richard Montgomery (UCSC) for many useful conversations and remarks.  Warm thanks to Corey Shanbrom (Cal State Sacramento) whose input greatly 
improved the overall exposition of the paper, Gary Kennedy (Ohio State), and Susan Colley (Oberlin).

\section{Preliminaries and Main Results}
Before presenting the results we want to give a brief summary of the $RVT$ coding system.  A more detailed discussion is presented in Subsection \ref{subsec:babymon}.
\

The $RVT$ coding system partitions the points at each level of the Monster Tower.  An $RVT$ code is a word in the letters $R$, $V$, $T_{i}$ for $i = 1,2$, and $L_{j}$ for $j = 1,2,3$ subject to certain
spelling rules.  The partial spelling rules were given in \cite{castro} and \cite{castro2} and state that the initial letter in any $RVT$ code must be the letter $R$, the letters $T_{i}$ and $L_{j}$ cannot
immediately follow the letter $R$, the first occurrence of the letter $L_{1}$ can only come after the letters $V$ and $T_{1}$, and the first appearance of the letters $T_{2}$ and $L_{j}$ for $j = 2,3$ cannot 
appear until the letter $L_{1}$ has been introduced.  This gives an incomplete list of spelling rules.  For example, we do not know what letters can come after 
$T_{2}$ or $L_{j}$ for $j = 2, 3$.

\begin{eg}
The codes $RRRR$, $RVT_{1}$, $RVL$, $RRVVT_{1}$, and $RVL_{1}T_{2}$ are allowable $RVT$ codes under our partial spelling rules.  The codes $RRT_{1}$, $RVRL_{2}$, and $RVT_{2}$ are not 
allowed, though. 
\end{eg}

We complete these spelling rules for the $\R ^{3}$-Monster Tower with the following theorems.  In Theorem \ref{thm:compspelling} the ``$:$'' denotes which letters can be placed after a
given letter.  For example, given the letter $R$ one can put either the letters $R$ or $V$ after it.       

\begin{thm}
\label{thm:compspelling}
The complete spelling rules for any $RVT$ code are as follows:
\begin{enumerate}
\item Any $RVT$ code string must begin with the letter $R$.
\item $R$ :  $R$ and $V$.
\item $V$ and $T \, (=T_{1})$ : $R, V, T,$ and $L$.
\item $L \, (= L_{1})$ and $L_{j}$ for $j = 1,2,3$ : $R, V, T_{i}$ for $i = 1,2$, and $L_{j}$ for $j = 1,2,3$.
\end{enumerate}
\end{thm}

In order to prove this we use Theorem \ref{thm:spelling} below as the base case of an induction argument showing that the spelling rules above hold for the first occurrence of
the distinguished directions $L_{j}$ for $j = 2,3$ and $T_{2}$.

\begin{thm}
\label{thm:spelling}
Let $\omega$ be an $RVT$ code of length $k$ where the $k$-th letter (the last letter in the code) is the first occurrence in the code of either the letters
$L_{j}$ for $j = 2,3$ or $T_{2}$.
\

\no Then the incidence relations between the various critical hyperplanes in $\Delta_{k}(p_{k})$ for $p_{k} \in \omega$ are determined and
summarized in Table \ref{tab:codes}.
\end{thm}

\begin{rem}[Theorem \ref{thm:spelling}]
We want to provide an example of what we mean by ``first occurrence'' in the statement of Theorem \ref{thm:spelling}.  The $RVT$ code $RVLT_{2}$ of length $4$ has the first occurrence of the letter 
$T_{2}$ as the last letter in the code.  Another example is the code $RVVLL_{3}$ of length $5$, where the first occurrence of the letter $L_{3}$ is the last letter of the code.   
\

Also, the critical hyperplanes will tell us which letters can come after the first occurrence of the letters $L_{j}$ for $j = 2,3$ and $T_{2}$.  By knowing which hyperplanes appear above these letter we are 
able to assertion that the letters $R$, $V$, $L_{j}$ for $j = 1,2,3$, and $T_{i}$ for $i = 1,2$ can come after the first occurrence of $L_{2}$ and $L_{3}$, and that
the letters $R$, $V$, $L_{3}$, and $T_{2}$ can be put after the first occurrence of the letter $T_{2}$.  This connection between the critical hyperplanes and $RVT$ code will be discussed in more detail in 
Subsection \ref{subsec:babymon}.  
\end{rem}

A {\it geometric distribution} hereafter denotes a linear subbundle of the tangent bundle with fibers of constant dimension.

\subsection{Prolongation}
Let the pair $(Z,\Delta)$ denote a manifold $Z$ of dimension $d$ equipped with a distribution $\Delta$ of rank $r$.  We denote by $\Pj (\Delta)$ the {\it projectivization} of $\Delta$, meaning
each vector space fiber $\Delta (p)$ for $p \in Z$ is projectivized resulting in fibers of dimension $(r-1)$.
As a manifold, $$\Pj (\Delta) \equiv Z^1,$$ has dimension $d + (r-1)$.

\begin{eg}
Take $Z = \mathbb{R}^3$, $\Delta = TR^3$ viewed as a rank 3 distribution. Then $Z^1$ is simply the trivial bundle $\mathbb{R}^3\times \mathbb{P}^2$, where the factor on the right denotes the
projective plane.
\end{eg}

Various geometric objects in $Z$ can be canonically prolonged (lifted) to the new manifold $Z^1$. In what follows prolongations of curves and transformations are essential.

The manifold $Z^1$ also comes equipped with a distribution $\Delta_1$ called the {\it Cartan prolongation of $\Delta$} (\cite{bryant}) which is defined as follows.
Let $\pi : Z^1 \rightarrow Z$ be the projection map $(p, \ell)\mapsto p$. Then
$$\Delta_1(p,\ell) = d\pi_{(p,\ell)}^{-1}(\ell),$$
i.e. {\it $\Delta_{1}(p, \ell)$ is the subspace of $T_{(p,\ell)}Z^1$ consisting of all tangents to curves in $Z$ that pass through $p$ with a velocity vector contained in $\ell$.}
It is easy to check using linear algebra that $\Delta_1$ is also a  distribution of rank $r$.

\begin{table}[t!]
\caption{Outline of Common Notation}
\begin{tabular}{|c|c|c|}
\hline
   Notation & Description   & Subsection \\
\hline
   $\mathcal{P}^{k}(n)$ & $k$-th level of the $\R^{n+1}$ Monster/Semple Tower & \ref{subsec:construct}   \\ \hline
   $\Delta_{k}$ & Rank $n$ distribution associated to $\mathcal{P}^{k}(n)$ & \ref{subsec:construct}    \\ \hline
   $\delta ^{j}_{i}(p)$ & Baby Monster through $p$ in $\Delta_{k}(p)$, $k = i + j$ & \ref{subsec:babymon} and \ref{sec:arrang} \\ \hline
   $\omega$ & The $RVT$ code of a point; string of $R$'s, $V$'s,$T$'s, and $L$'s & \ref{subsec:construct},
      \ref{subsec:rc}, and \ref{subsec:babymon} \\ \hline
\end{tabular}
\label{tab:notation}
\end{table}

We note that the word prolongation will always be synonymous with Cartan prolongation.

\subsection{The Monster Tower and $RVT$ Coding.}
\label{subsec:construct}
We start with $\R^{n+1}$ as our base manifold $Z$ and take $\Delta_{0} = T \R^{n+1}$.  Prolonging $\Delta_{0}$ we get
$\mathcal{P}^{1}(n) = \Pj(\Delta_{0})$ equipped with the distribution $\Delta_{1}$ of rank $n+1$.  By
iterating this process we end up with the manifold $\mathcal{P}^{k}(n)$ which is endowed with the rank $n+1$
distribution $\Delta_{k} = (\Delta_{k-1})^{1}$ and fibered over $\mathcal{P}^{k-1}(n)$.
In this paper we will be studying the case $n=2$.

\begin{defn}
The \textit{$\R^{n+1}$-Monster Tower} is a sequence of manifolds with distributions, $(\mathcal{P}^{k}, \Delta_{k})$,
together with fibrations $$\cdots \rightarrow \mathcal{P}^{k}(n) \rightarrow \mathcal{P}^{k-1}(n) \rightarrow \cdots \rightarrow
\mathcal{P}^{1}(n) \rightarrow \mathcal{P}^{0}(n) = \mathbb{R}^{n+1}$$
and we write $\pi_{k,i}: \mathcal{P}^{k}(n) \rightarrow \mathcal{P}^{i}(n)$ for the respective bundle projections.
\end{defn}

\begin{table}[t!]
\caption{Critical Hyperplanes Incidence Relations.}
\begin{tabular}{|c|c|c|}
\hline
$RVT$ Code &      $T_{1}$ & $T_{2}$ \\
\hline
$R \omega L_{1}$            &        $\delta ^{1}_{k-1}(p_{k})$       &  $\delta ^{2}_{k-2}(p_{k})$  \\
\hline
$R \omega L_{2}$:          &             &                                                                                                       \\  \hdashline

   $R \omega' VT^{m} LL_{2} \, \, \text{for} \, \, m \geq 0$           &        $\delta ^{2}_{m+2+r}(p_{m+r+4})$      &  $\delta ^{m+3}_{1 + r}(p_{m+r+4})$             \\
   $R \omega' LLL_{2}$                                                                 &       $\delta ^{2}_{k-2} (p_{k})$      & $\delta ^{3}_{k-3}(p_{k})$            \\
\hline
$R \omega L_{3}$:            &            &                                                                                                                                                                     \\   \hdashline

   $R \omega' VT^{m} L L_{3} \, \, \text{for} \, \, m \geq 0$           &        $\delta ^{1}_{k-1}(p_{k})$ &   $\delta ^{m + 3}_{r+1}(p_{r+m+4})$  \\
   $R \omega' LL L_{3}$                                                                &        $\delta ^{1}_{k-1}(p_{k})$ &   $\delta ^{3}_{k-3}(p_{k})$   \\
   $R \omega' V T^{m}L T_{2} L_{3} \, \, \text{for} \, \, m \geq 0$  &  $\delta ^{1}_{k-1}(p_{k})$ & $\delta ^{m + 4}_{r + 1}(p_{m+r+5})$ \\
   $R \omega' LLT_{2}L_{3} $ & $\delta ^{1}_{k-1}(p_{k})$  &    $\delta ^{4}_{k-4}(p_{k})$                                                                           \\
   $R \omega' V T^{m}L L_{2} T_{2} L_{3} \, \, \text{for} \, \, m \geq 0$  &  $\delta ^{1}_{k-1}(p_{k})$   &    $\delta ^{4}_{m+r+2}(p_{m+r+6})$   \\
   $R \omega' V T^{m} L ^{s}L_{2}T_{2} L_{3} \, \, \text{for} \, \, m \geq 0, s \geq 2$  &   $\delta ^{1}_{k-1}(p_{k})$   &  $\delta ^{4}_{m+r+2}(p_{m+r+6})$    \\
\hline
$R \omega T_{1}$:           &              &                                                                                                                                                                     \\    \hdashline
    $R \omega' L T_{1}$                                                                 &            $\delta ^{2}_{k-2}(p_{k})$  &   None                                        \\
    $R \omega' L_{1} T_{1}$                                                          &            $\delta ^{2}_{k-2}(p_{k})$ &    None                                        \\
    $R \omega' L_{2} T_{1}$                                                          &            $\delta ^{3}_{k-3}(p_{k+3})$  &    None                                         \\
    $R \omega' L_{3} T_{1}$                                                          &            $\delta ^{2}_{k-2}(p_{k})$       &    None                                         \\
\hline
$R \omega T_{2}$:           &              &                                                                                                                                                                        \\    \hdashline
  $R \omega' V T^{m} L T_{2}$ \, \, \text{for} \, \,  $m \geq 0$  &          None                                    & $\delta ^{m+3}_{1+r}(p_{m+ r +4 })$      \\
  $R \omega' L L T_{2}$                                                           &          None                                    & $\delta ^{3}_{k-3}(p_{k})$                  \\
  $R \omega' V T^{m} L L_{2}T_{2} \, \, \text{for} \, \, m \geq 0$ &    None                                          & $\delta ^{m+4}_{r+1}(p_{m+r+5})$            \\
  $R \omega' V T^{m}L^{s}L_{2}T_{2} \, \, \text{for} \, \, m \geq 0, s \geq 2$  &  None                      & $\delta ^{4}_{r+m+s}(p_{m+s+r+4})$        \\
  $R \omega' VT^{m}LL_{3}T_{2} \, \, \text{for} \, \, m \geq 0$  &   None                      & $\delta ^{m+4}_{r+1}(p_{m+r+5})$                            \\
  $R \omega' V T^{m} L^{s}L_{3} T_{2} \, \, \text{for} \, \, m \geq 0, s \geq 2$ & None                & $\delta ^{3}_{m+s+r}(p_{m+s+r+3})$            \\
\hline
\end{tabular}
\label{tab:codes}
\end{table}

\begin{rem}[Notation from Table \ref{tab:codes}]
We assume that each of the above $RVT$ codes is of total length $k$.  At the same time, each of the intermediate and undetermined $RVT$ blocks, denoted by $\omega'$, are assumed to be of total 
length $r$.
\

It is worth noting that the $L_{2}$ and $L_{3}$ directions will have the same critical hyperplane configuration as $L$ points, but the level at which the Baby Monsters are born from which will produce 
these critical hyperplanes will be different than in the $L$ case.  In short, the source of the critical hyperplanes over the $L_{j}$'s depends on which letters come before these distinguished directions.
\

When we prove Theorem \ref{thm:compspelling} we can restrict our attention to looking at the cases of $L_{j}$ for $j = 2,3$ and $T_{2}$, since from Definition \ref{defn:Tpoints} (below) we can say 
$T = T_{1}$ and $L = L_{1}$, which from previous work, we know what the spelling rules are for these letters.  Despite this fact, we will still list where the critical hyperplane come from for these
cases in Table \ref{tab:codes}.
\

The letter $T_{2}$ gives us a new hyperplane configuration that was previously unknown and is illustrated in Figure \ref{fig:t2-planes}.
\end{rem}

In proving Theorem \ref{thm:spelling}, we will present a method to determine what critical hyperplanes and incidence relations will occur over any point in the $\R ^{3}$-Monster
Tower.  This method can be applied to the $\R^{n}$ Monster Tower as well.  In addition, this algorithm will tell us which fibers the Baby Monsters are born from.
\

We also want to point out that in Subsection \ref{subsec:rc} we will provide the basic definitions needed to understand the $RVT$ coding system.
Since we will be working exclusively with the $n=2$ Monster Tower in this paper, we will just write
$\mathcal{P}^{k}$ for $\mathcal{P}^{k}(2)$.






\subsection{Baby Monsters.}
\label{subsec:babymon}
One can apply prolongation to any analytic $m$-dimensional manifold $F$ in place of $\R^{m}$.  Start out with
$\mathcal{P}^{0}(F) = F$ and take $\Delta^{F}_{0} = TF$.  Then the prolongation of the pair
$(F, \Delta^{F}_{0})$ is $\mathcal{P}^{1}(F) = \Pj TF$ equipped with the rank $m$ distribution
$\Delta^{F}_{1} \equiv (\Delta^{F}_{0})^{1}$.  By iterating this process $k$ times we end up with new the pair
$(\mathcal{P}^{k}(F), \Delta^{F}_{k})$, which is analytically diffeomorphic to
$(\mathcal{P}^{k}(m-1), \Delta_{k})$ (\cite{castro}).
\

Now, apply this process to the fiber $F_{i}(p_{i}) = \pi^{-1}_{i, i-1}(p_{i-1}) \subset \mathcal{P}^{i}(m-1)$ through the point
$p_{i}$ at level $i$. The fiber is an $(m-1)$-dimensional integral submanifold for $\Delta_{i}$.  Prolonging, we see that
$\mathcal{P}^{1}(F_{i}(p_{i})) \subset \mathcal{P}^{i + 1}(m-1)$, and $\mathcal{P}^{1}(F_{i}(p_{i}))$ has the associated distribution
$\delta^{1}_{i} \equiv \Delta^{F_{i}(p_{i})}_{1}$; that is,
$$\delta^{1}_{i}(q) = \Delta_{i + 1}(q) \cap T_{q}(\mathcal{P}^{1}(F_{i}(p_{i}))) $$
which is a hyperplane within $\Delta_{i + 1}(q)$, for $q \in \mathcal{P}^{1}(F_{i}(p_{i}))$.  When this prolongation process is
iterated, we end up with the submanifolds
$$\mathcal{P}^{j}(F_{i}(p_{i})) \subset \mathcal{P}^{i + j}(m-1) $$
with the hyperplane subdistribution $\delta^{j}_{i}(q) \subset \Delta_{i + j}(q)$ for $q \in \mathcal{P}^{j}(F_{i}(p_{i}))$.

\begin{defn}   A {\it Baby Monster} born at level $i$, is a sub-tower $(\mathcal{P}^{j}(F_{i}(p_{i})), \delta^{j}_{i})$,
for $j \geq 0$ within the ambient Monster Tower.  If $q \in \mathcal{P}^{j}(F_{i}(p_{i}))$ then we will say that a Baby Monster born
at level $i$ passes through $q$ and that $\delta^{j}_{i}(q)$ is a \textit{critical hyperplane} passing through $q$, which was born at level $i$.
\end{defn}

\begin{rem}
In this paper $n = 3$, where our base is $\R^{3}$, which results in our the distribution $\Delta_{i+j}$ being of rank $3$ and that all of our critical hyperplanes $\delta ^{j}_{i}$ being $2$-planes within the 
distribution $\Delta_{i+j}$.  As a result, we will replace the term ``hyperplane'' by just``plane''  for the rest of the paper.
\end{rem}

\begin{defn} The \textit{vertical plane} $V_k (q)$, often written as just $V(q)$, is the critical plane $\delta^{0}_{k} (q) = T_{q} (F_{k}(p_{k}))$.  We note that it is always one of the critical planes 
passing through $q$.
\end{defn}

\begin{figure} 
  \centering
   \subfloat[Above a regular point.]{\label{fig:one-plane}\includegraphics[clip = true, width=0.36\textwidth]{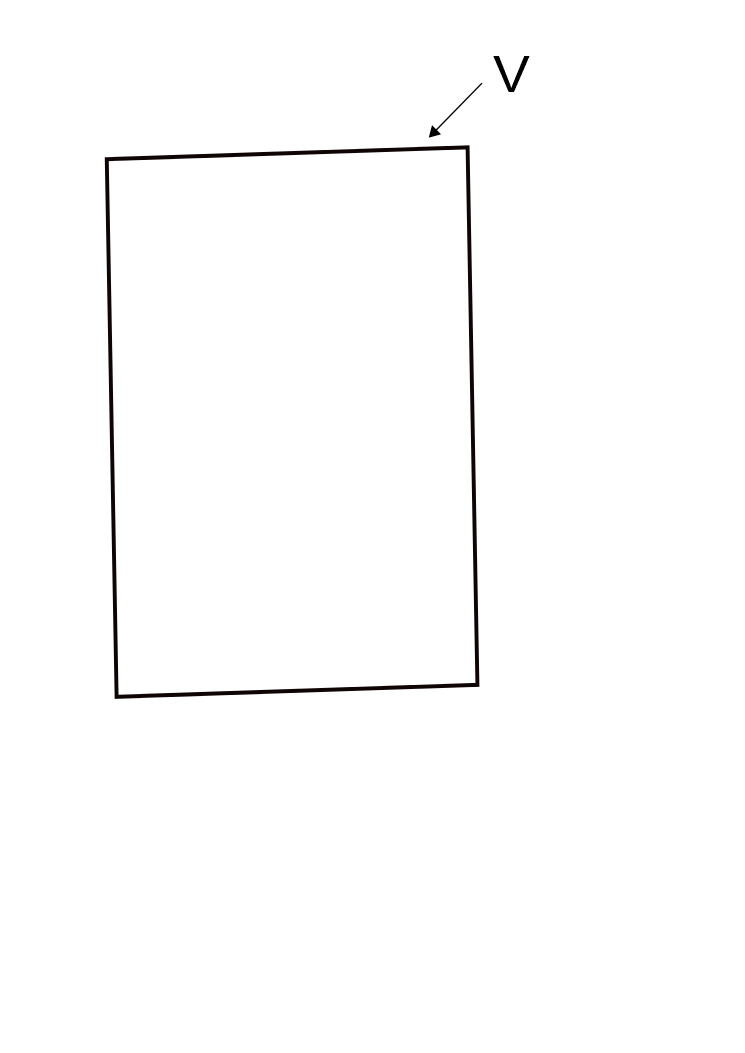}}
   \subfloat[Above a vertical or tangency point.]{\label{fig:two-planes}\includegraphics[width=.39\textwidth]{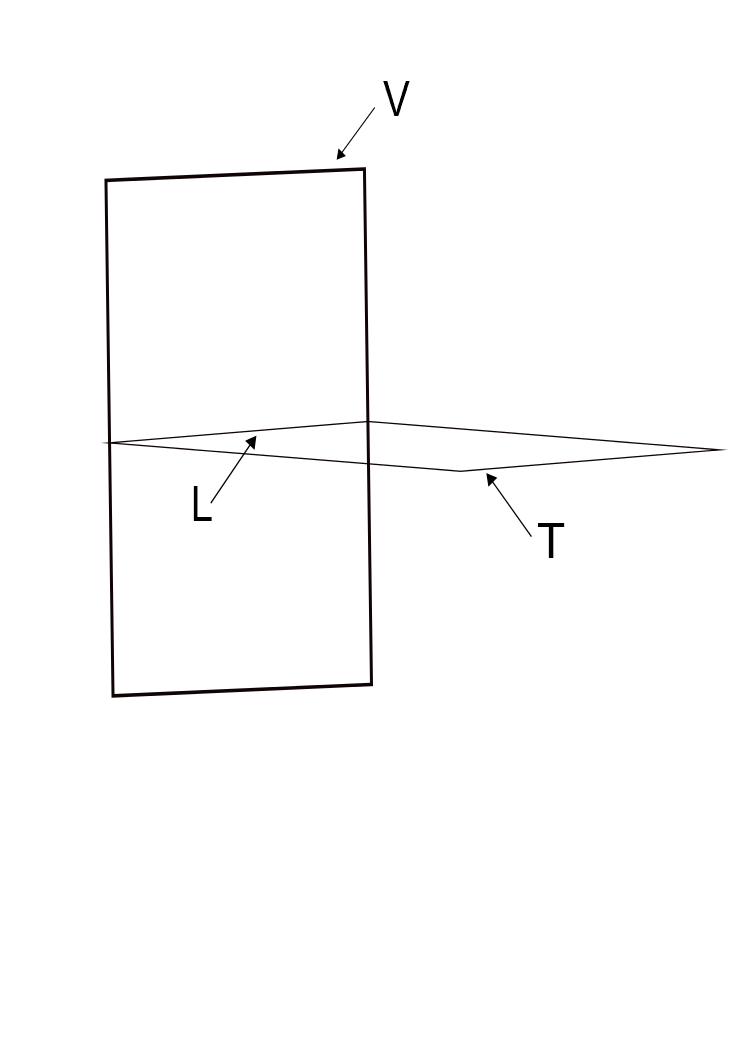}}
   \caption{Arrangement of critical planes over the letters $R$, $V$, and $T$.}
   \label{fig:arrangement}
\end{figure}

\subsection{$RC$ coding of points.}
\label{subsec:rc}

\begin{defn}
A direction $\ell \subset \Delta_{k}(p_{k})$, $k \geq 1$ is called a {\it critical direction} if there exists an
immersed curve at level $k$ that is tangent to the direction $\ell$, and whose projection to level zero, meaning
the base manifold, is a constant curve.
If no such curve exists,
then we call $\ell$ a {\it regular direction.}  Note that while $\ell$ is technically a line we will by an
abuse of terminology refer to it as a direction.
\end{defn}

\begin{defn}
A point $p_{k} \in \mathcal{P}^{k}$, where $p_{k} = (p_{k-1}, \ell)$ is called a \textit{regular or critical point} if the line $\ell$ is a regular direction or a critical direction.
\end{defn}

\begin{defn}
For $p_{k} \in \mathcal{P}^{k}$, $k \geq 1$ and $p_{i} = \pi_{k,i}(p_{k})$, we write
$\omega_{i}(p_{k}) = R$ if $p_{i}$ is a regular point and $\omega_{i}(p_{k}) = C$ if $p_{i}$ is a critical point.
Then the word $\omega(p_{k}) = \omega_{1}(p_{k}) \cdots \omega_{k}(p_{k})$ is called the $RC$ code for the point $p_{k}$.  The number of letters within the $RC$ code for $p_{k}$ equals the level of
the tower that the
point lives in.  Note that $\omega_{1}(p_{k})$ is always equal to $R$, see \cite{castro2}.
\end{defn}

\no In the following section we will show
that there is more than one kind of critical direction that can appear within the distribution $\Delta_{k}$.

\subsection{Arrangements of critical planes for $n = 2$.}\label{sec:arrang}
\begin{figure} 
  \centering
  \subfloat[Above an $L_{j}$ for $j = 1,2,3$ point.]{\label{fig:three-planes}\includegraphics[width=0.45\textwidth]{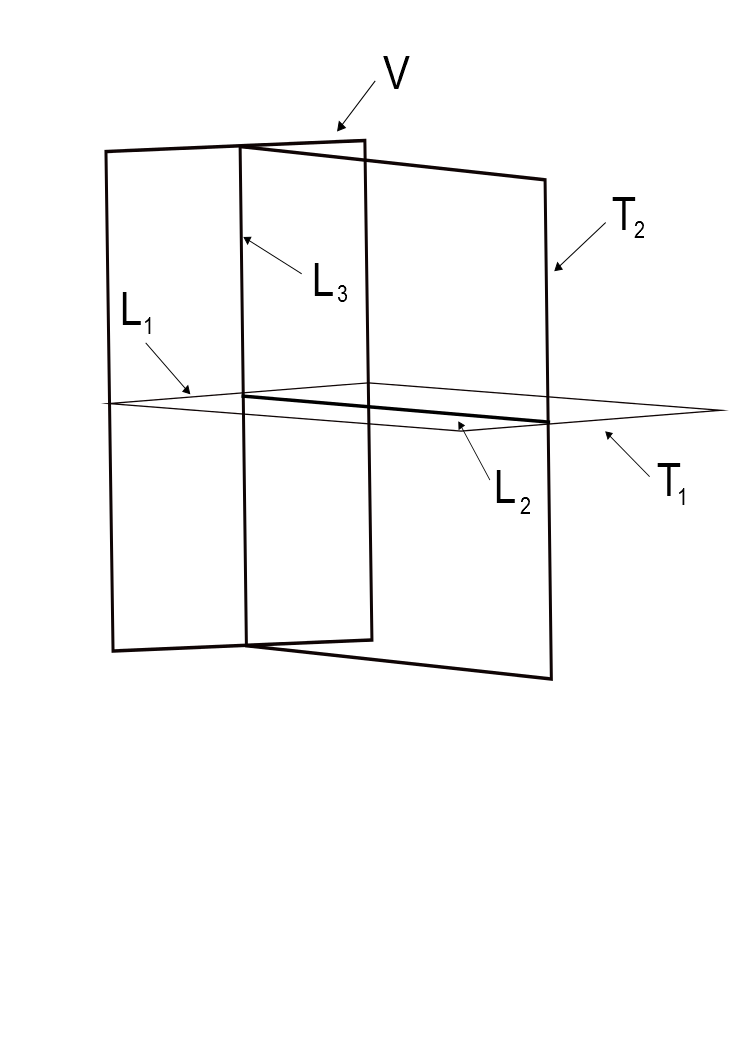}}
  \subfloat[Above a $T_{2}$ point.]{\label{fig:t2-planes}\includegraphics[width=0.45\textwidth]{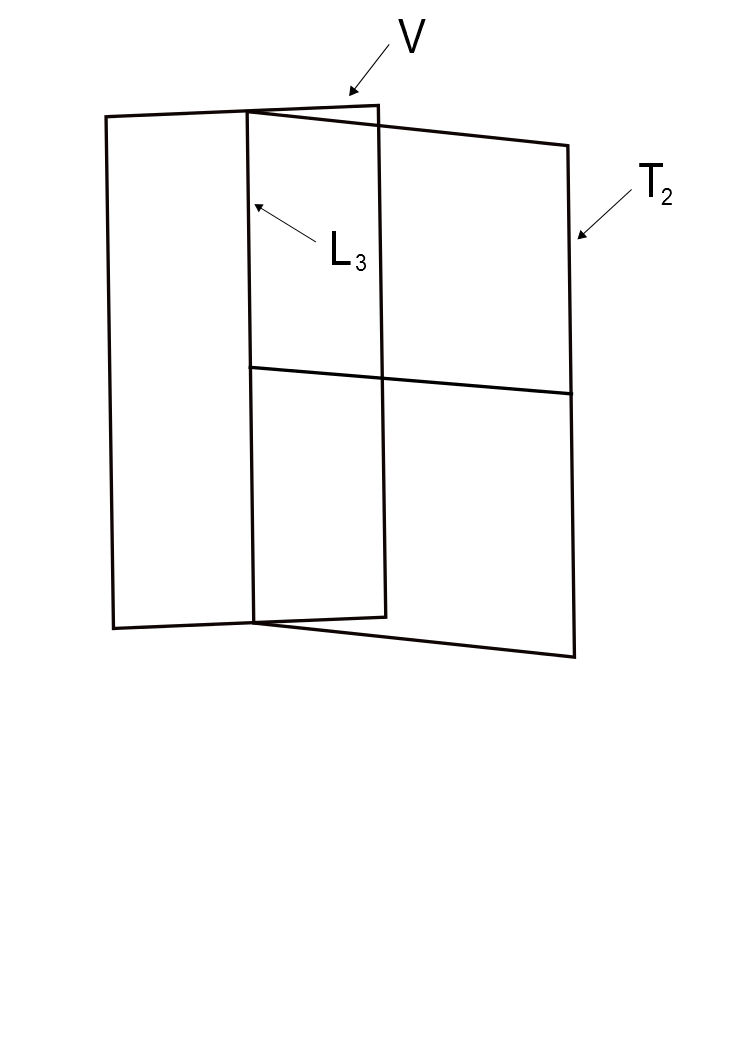}}
  \caption{Arrangement of critical planes over $L_{j}$ for $j = 1,2,3$ and $T_{2}$.}
  \label{fig:arrangement2}
\end{figure}

Over any point $p_{k}$, at the $k$-th level of the Monster Tower, there is a total of four possible plane configurations
for $\Delta_{k}$.  These three configurations are shown in Figures \ref{fig:one-plane}, \ref{fig:two-planes}, \ref{fig:three-planes}, and \ref{fig:t2-planes}.

Figure $\ref{fig:one-plane}$ is the picture
for $\Delta_{k}(p_{k})$ when the $k$-th letter in the $RVT$ code for $p_{k}$ is the letter $R$.  This means that the
vertical plane, labeled with a $V$, is the only critical plane sitting inside of $\Delta_{k}(p_{k})$.  Figure $\ref{fig:two-planes}$ is the
picture for $\Delta_{k}(p_{k})$ when the $k$-th letter in the $RVT$ code is either the letter $V$ or the letter $T$.  This gives a total of two critical planes sitting inside of $\Delta_{k}(p_{k})$ and one
distinguished critical direction: one is the vertical plane and the other is the
tangency plane, labeled by the letter $T$.  The intersection of vertical and tangency plane gives a
distinguished critical direction, which is labeled by the letter $L$.  Now, Figure $\ref{fig:three-planes}$ describes the picture for $\Delta_{k}(p_{k})$ when the $k$-th letter
in the $RVT$ code of $p_{k}$ is the letter $L$.  Figure $\ref{fig:three-planes}$ depicts this situation where there is now a total of three
critical planes: one is the vertical plane, and two tangency planes, labeled as $T_{1}$ and $T_{2}$.  Now, because of
the presence of these three critical planes we need to refine our notion of an $L$ direction and add two more distinct $L$ directions.
These three directions are labeled as $L_{1}$, $L_{2}$, and $L_{3}$.  We will prove in this paper that this configuration will also persist for $L_{j}$ for $j = 2,3$, as well.  The last figure given in
Figure \ref{fig:t2-planes} is the critical plane configuration over the letter $T_{2}$.  It is worth pointing out that this is an unexpected result, at least to the authors of the paper,
and gives us a new incidence relation that was previously unknown.

\

With the above pictures in mind, we can now refine our $RC$ coding and define the $RVT$ code for points within the Monster Tower.
Take $p_{k} \in \mathcal{P}^{k}$ and if $\omega_{i}(p_{k}) = C$ then we look at the point $p_{i} = \pi_{k,i}(p_{k})$, where
$p_{i} = (p_{i-1}, \ell_{i-1})$.  Then depending on which critical plane, or distinguished direction, contains $\ell_{i-1}$, we replace the letter $C$ by one of
the letters $V$, $T_{i}$ for $i = 1,2$, or $L_{j}$ for $j = 1,2,3$.
We need to define the two distinguished tangency planes and the three
distinguished directions that arise from the intersection of the three critical planes when the two tangency planes are present in $\Delta_{k}$.

\subsection{$KR$-Coordinates}
\label{subsec:kr}
For our work within the Monster Tower one needs to work with a suitable coordinate system.  In Section $4.4$ of \cite{castro2} we detailed the basics of \textit{Kumpera-Rubin coordinates}, or
$KR$-coordinates for short, by determining the coordinates for the $RVT$ class $RVL$.  We will briefly summarize this construction to help the reader understand some of the basic properties
of this coordinate system.  
\

We begin with the pair $(\mathcal{P}^{0}, \Delta_{0}) = (\R ^{3}, T \R ^{3})$ with the coframe $\{ dx, dy, dz \}$ for $T \R ^{3}$.  We center our chart at $(p_{0}, \ell_{0})$ where 
$p_{0} = (0,0,0)$ and $\ell_{0}$ is a direction in $T_{p_{0}} \R ^{3}$ such that $dx \rest{\ell_{0}} \neq 0$.  This allows us to introduce fiber affine coordinates
$[dx : dy : dz] = [1 : \frac{dy}{dx} : \frac{dz}{dx} ]$, where
$$u = \frac{dy}{dz} \, \, \, \, v = \frac{dz}{dx} ,$$
\no and results in 
$$\{ dy - udx = 0, \, \, dz - vdx = 0  \} \, = \, \Delta_{1}.$$
\no Now take the point $p_{1} = (p_{0}, \ell _{0})$ in the first level of the tower $\mathcal{P}^{1}$ and look at the line $\ell_{1} \subset \Delta_{1}(p_{1}')$ for $p_{1}'$.  Recall that the vertical
plane is given by $V = \delta ^{0}_{1}(p_{1}) = \mbox{span} \{ \frac{\pa}{\pa u}, \frac{\pa}{\pa v} \}$.  Let us suppose $\ell_{1} = \mbox{span} \{ \frac{pa}{\pa u} \}$ so that it is in the vertical plane.
Then near $\ell_{1}$ we have $[dx : du : dv] = [\frac{dx}{du} : 1 : \frac{dv}{du}]$ to give the affine coordinates
$$u_{2} = \frac{dx}{du} \, \, \, \, v_{2} = \frac{dv}{du} ,$$
\no with 
\begin{align*}
\Delta_2 = \{&dy - u dx = 0, dz - v dx = 0,\\
&dx - u_2 du = 0, dv - v_2 du = 0\} .
\end{align*} 

Then the rank $3$ distribution $\Delta_{2}$ is coframed by $[du : du_{2} : dv_{2}]$ with the vertical plane given by $du = 0$ and the tangency plane $T$ given by
$du_{2} = 0$.  The point $p_{3} = (p_{2}, \ell)$ with $\ell$ being an $L$ direction means that both $du\rest{\ell} = 0$ and $du_{2}\rest{\ell} = 0$.  As a result, the only choice for
local coordinates near $p_{3}$ is given by $[\frac{du}{dv_{2}}: \frac{du_{2}}{dv_{2}}: 1]$ to give the fiber coordinates  
$$u_{3} = \frac{du}{dv_{2}}  \, \, \, \, v_{3} = \frac{du_{2}}{dv_{2}} ,$$
\no where the distribution $\Delta_{3}$ is described by
\begin{align*}
\Delta _{3} = \{&dy - u dx = 0, dz - v dx = 0,\\
               &dx - u_2 du = 0, dv - v_2 du = 0, \\
               &du - u_{3}dv_{2} = 0, du_{2} - v_{3}dv_{2} = 0 \}.
\end{align*}
\no From here one can show that $T_{1} = \delta ^{1}_{2}(p_{3})$ and $T_{2} = \delta ^{2}_{1}(p_{3})$ in $\Delta_{3}(p_{3})$ with $p_{3}$ given in $KR$-coordinates by
$p_{3} =(x, y, z, u, v, u_{2}, v_{2}, u_{3}, v_{3}) = (0, 0, 0, 0, 0, 0, 0, 0, 0)$.

\begin{defn}[Definitions of the $T_{1}$ and $T_{2}$ critical planes]
\label{defn:Tpoints}
At any given level of the Monster Tower, say the $k$-th level, we can write the fiber space at the $k$-th level in the form $F_{k}(p_{k}) = \pi ^{-1}_{k,k-1}(p_{k-1})$ and is written
in $KR$-coordinates, see \cite{castro} for details, as $F_{k}(p_{k}) = (p_{k-1}, u_{k}, v_{k})$, where $p_{k-1}$ is a fixed point at the $(k-1)$-st level.  We can represent the vertical space $V_{k}$  
in a neighborhood of the point $p_{k}$ as $\delta ^{0}_{k}(p_{k}) = span \{ \frac{\pa}{\pa u_{k}}, \frac{\pa}{\pa v_{k}} \}$.  $T_{1}$ is the critical plane that intersects $span \{ \frac{\pa}
{\pa v_{k}} \}$.  In local $KR$-coordinates for $\Delta_{k}$ it is given by $[df_{k}: du_{k} : dv_{k}] = [a : 0 : b]$ for $a,b \in \R$, with $a \neq 0$. \footnote{The covector $df_{k}$ is known as the 
\textit{uniformizing coordinate}, see \cite{castro} for more details.}
\

The $T_{2}$ directions are characterized as the directions which do not intersect the $\frac{\pa}{\pa v_{k}}$ component of the vertical space 
$V_{k}$.  In $KR$-coordinates we have the $T_{2}$ critical plane is characterized by $[df_{k} : du_{k} : dv_{k}] = [a : b : 0]$ for $a,b \in \R$ and $a \neq 0$.
\end{defn}

\begin{defn}[Definition of the $L_{j}$ directions for $j = 1,2,3$]
The distinguished line $L_{1}$ is given by $V \cap T_{1}$, $L_{2}$ by $T_{1} \cap T_{2}$, and $L_{3}$ by $V \cap T_{2}$.
\end{defn}

\begin{eg}[Examples of $RVT$ codes]
The following are examples of $RVT$ codes: $R \cdots R$, $RVVT$, $RVLT_{2}R$, and $RVLL_{2}$.
The code $RTL$ is not allowed because the letter $T$ is preceded by the letter $R$ and the code $RLT_{2}$ is
not allowed because the letter $L$ comes immediately after the letter $R$.
\end{eg}

\begin{figure}
 \label{fig:exmpRVL}
   \centering
   \includegraphics[width=0.57\textwidth]{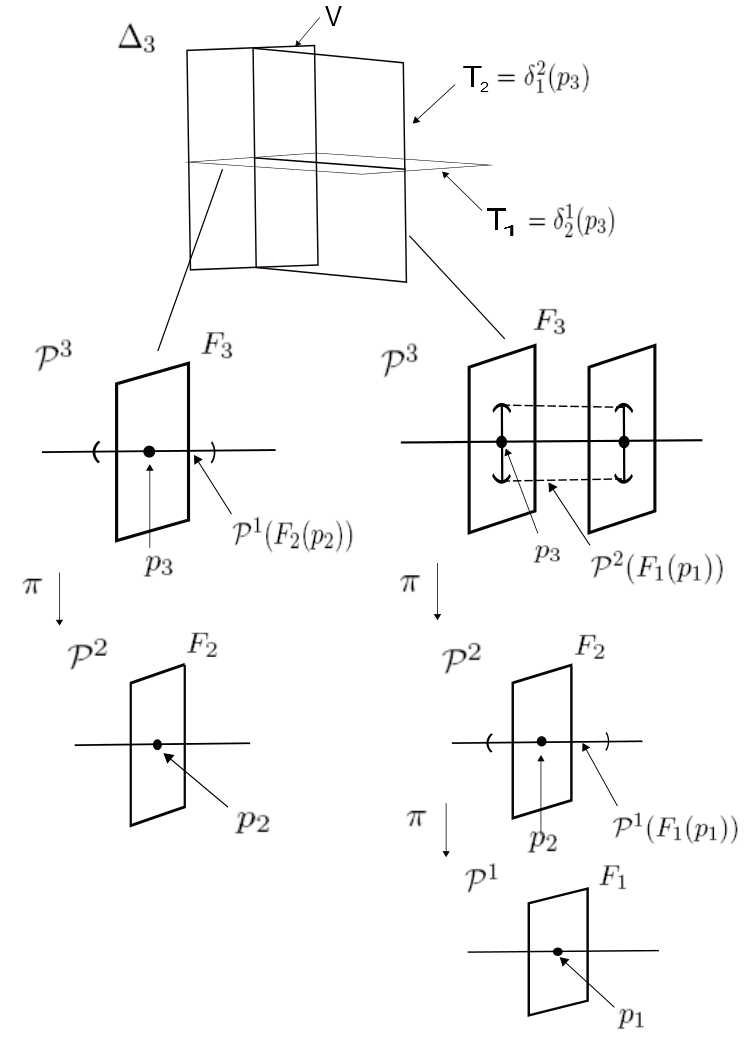}
   \caption{Critical plane configuration over $p_{3} \in RVL$.}
   \label{fig:overl}
\end{figure}

In the subsequent example we will redo the calculation for the critical planes over the $L$ direction that was first presented in \cite{castro2}.
This time though we will be working with a different method which tells us where the various critical planes come from in the distribution $\Delta_{k}$ over an
$L$ point.  This will help us avoid guessing which levels these critical planes originate from.  We will also be using this method in order to prove Theorem \ref{thm:spelling}.  
Figure \ref{fig:overl} illustrates which fibers these critical planes come from and our incidence relations.

\begin{eg}[Calculation of Incidence Relations for the Class $RVL$]
\label{eg:levelthree}
We begin with the local $KR$-coordinates for a point $p_{3}$ in the third level of the $\R ^{3}$- Monster Tower, which is in the $RVT$ class $RVL$.  From the work done in Example $4.4$ of 
\cite{castro2} we have that $\Delta_{3}$ in a neighborhood of $p_{3} \in RVL$ is coframed by
$[dv_{2} : du_{3} : dv_{3}]$.  The vertical plane is characterized by $[0 : a : b]$ for $a,b \in \R$ with $a \neq 0$ in these coordinates.  If there are any other critical planes then they are either of 
the form $[a : 0 : b]$ with $a \neq 0$ to give the critical plane $T_{1}$ or of the form $[a : b : 0]$ with $a \neq 0$ to give the $T_{2}$ critical plane.  We will start with planes of 
this form and work backwards until we determine which fiber they must have come from.  If we are unable to do this or arrive at a contradiction along the way (such as showing that there is a tangency 
plane above a regular point), then it implies that no such critical plane could have existed at that level.
\

We start with $[dv_{2} : du_{3} : dv_{3}] = [a : 0 : b]$ and work backwards.  From looking at these projective coordinates it tells us that the coordinates $v_{2}$ and $v_{3}$ can't both be zero in a
neighborhood of $p_{3}$ for the Baby Monster that gives rise to this critical plane, but $u_{3}$ must be identically zero.  Then, one level down we have that $\Delta_{2}$ is coframed in a
neighborhood of $p_{2} \in RV$ by $[du : du_{2} : dv_{2}]$ and since we move in an $L$ direction the coordinates on the third level are $u_{3} = \frac{du}{dv_{2}}$ and $v_{3} = \frac{du_{2}}{dv_{2}}$.
This implies that the $u$ coordinate must be identically zero on the Baby Monster and while the $v_{3}$ coordinate will be zero at $p_{3}$ it cannot be identically zero on the manifold.  As a result,
our Baby Monster one level down is given by $[du : du_{2} : dv_{2}] = [0 : a : b]$.  This tells us we can stop and do not need to go any farther down in the tower because our critical plane comes
from the fiber space $F_{2}(p_{2}) = (p_{1}, u_{2}, v_{2})$, implying that $T_{1}$ comes from the Baby Monster $\delta ^{1}_{2}$.
\

The second tangency plane is given by $[dv_{2} : du_{3} : dv_{3}] = [a : b : 0]$.  Again, these projective coordinates tell us that the coordinates $v_{2} \not \equiv 0$ and $u_{3} \not \equiv 0$ in a 
neighborhood of $p_{3}$ for our Baby Monster which gives rise to the $T_{2}$ tangency plane, but $v_{3} \equiv 0$.
One level down, in a neighborhood of $p_{2}$, we have $\Delta_{2}$ is coframed by $[du: du_{2} : dv_{2}]$, since we moved in an $L$ direction one level up and our plane is given by 
$[du : du_{2} : dv_{2}] = [a : 0 : b]$.  This implies that the $u_{2} \equiv 0$ along the Baby Monster.  Going down one more level we have $\Delta_{1}$ 
in a neighborhood of $p_{1}$ is coframed by $[dx : du : dv]$.  Since we moved in a $V$ direction we end up with $[dx : du : dv] = [\frac{dx}{du} : 1 : \frac{dv}{du}]$ to give the coordinates.
$u_{2} = \frac{dx}{du}$ and $v_{2} = \frac{dv}{du}$, which results in $[dx : du : dv] = [0 : a : b]$ along our Baby Monster.  We can stop at this point and see that we do not need to go any farther down
in the tower because our Baby Monster will arise from the fiber space $F_{1}(p_{1}) = (p_{0}, u, v)$, indicating that $T_{2}$ comes from the Baby Monster $\delta ^{2}_{1}$, which again is in agreement 
with what we expect from Example $4.4$ in \cite{castro2}.
\end{eg}

\section{The Critical Plane Method}
The above example helps to illustrate our algorithm for determining the various critical planes.  Here is the overall idea of our algorithm:  Take an $RVT$ code of length $k$ called $\omega$.
Say we are interested in understanding which critical letters we can add to the end of $\omega$, and, at the same time, understanding where the critical planes originate from within the Monster
Tower.  Take a point $p_{k} \in \omega$ and determine what the local $KR$-coordinates are that describe the distribution $\Delta_{k}(p_{k})$.  From there we start with how the
various critical planes are presented within the coframing for the distribution $\Delta_{k}(p_{k})$, based upon Definition \ref{defn:Tpoints}.  Then we work backwards looking at the way these
critical planes are represented in the coframing for each $\Delta_{i}(p_{i})$ for $0 \leq i \leq k-1$.  If we end up with two fiber coordinates being nontrivial for our critical plane at some
intermediate level, meaning $u_{j}, v_{j}$ for some $1 \leq j \leq k-1$ will both be nonvanishing for our critical plane, it tells us that our critical plane within $\Delta_{k}(p_{k})$ will originate from 
the fiber space $F_{j}(p_{j})$.  If we can't find such a fiber space or if this critical plane is  projected to a lower level and gives rise to a critical plane in $\Delta_{i}(p_{i})$ for $1 \leq i \leq k-1$ 
which can't appear in the distribution based upon our spelling rules, then it will imply that this critical plane can not appear within the distribution $\Delta_{k}(p_{k})$.
\

While the above does not give an algorithm in the strict sense, it is a method for how one can determine if these critical planes exist or not within our distribution.  
We believe that it will also be much more illuminating to present the mechanics of how the above outline works through various examples such as the one above for the case of $RVL$ and the ones that 
follow below.

\subsection{The Critical Planes over the $T_{2}$ Plane}
\label{subsec:exampleT2}
This section will be building upon the calculations done in Example \ref{eg:levelthree}.  We show that the critical plane configuration over any point strictly in the code $RVLT_{2}$ in any critical plane is 
the vertical plane and the $T_{2}$ plane shown in Figure $5$. 
\

The distribution $\Delta_{3}$ 
in a neighborhood of $p_{3} \in RVL$ is coframed by $[dv_{2}: du_{3}: dv_{3}]$ and for $\ell \subset T_{2}$ we have the $dv_{3} = 0$ and must have $dv_{2} \neq 0$ $du_{3} \neq 0$.  This results in
fiber coordinates $u_{4} = \frac{du_{3}}{dv_{2}}$ and $v_{4} = \frac{dv_{3}}{dv_{2}}$ with $u_{4} \neq 0$ and $v_{4} = 0$.  For a direction $\ell$ to be completely in the $T_{2}$ plane it needs to be a linear 
combination of both the directions $\frac{\pa}{\pa v_{2}}$ and $\frac{\pa}{\pa u_{3}}$, i.e. in the $span \{ \frac{\pa}{\pa v_{2}}, \frac{\pa}{\pa u_{3}} \}$.  The distribution in a neighborhood of $p_{4} \in 
RVLT_{2}$ is given by
\begin{align*}
\Delta_{4} = \{&dy - u dx = 0, dz - v dx = 0,\\
                         &dx - u_2 du = 0, dv - v_2 du = 0 \\
                         &du - u_{3}dv_{2} = 0, du_{2} - v_{3}dv_{2} = 0 \\
                         &du_{3} - u_{4}dv_{2} = 0, dv_{3} - v_{3}dv_{2} = 0 \} \subset T\mathcal{P}^{4}.
\end{align*}
\

Our next step is to determine what the critical planes are below the point $p_{4}$ in the distribution $\Delta_{4}(p_{4})$.  In short, one can look at the prolongation of the fibers $F_{3}(p_{3})$,
$F_{2}(p_{2})$, and $F_{1}(p_{1})$ to see what the critical planes are and if they will appear in the distribution over the point $p_{4}$.  However, for a more general $RVT$ code it is difficult and 
somewhat cumbersome to prolong each fiber space.  Even in this case one would have to prolong $3$ different fiber spaces in order to determine where the possible critical planes would come 
from.  We will show, using the algorithm presented in Example \ref{eg:levelthree}, that the only critical plane, other than the vertical one, over the point $p_{4} \in RVLT_{2}$ is the $T_{2}$ critical 
plane that comes from the fiber $F_{1}(p_{1})$.

\subsubsection{Showing that $T_{2} = \delta ^{3}_{1}(p_{4})$}
We begin with the coframing for $\Delta_{4}$ in a neighborhood of the point $p_{4}$ given in $KR$-coordinates by $[dv_{2} : du_{4} : dv_{4}]$.  If the $T_{2}$ plane is to exist it must be of the form
$[dv_{2} : du_{4} : dv_{4}] = [a : b : 0]$ for $a,b \in \R- {0}$.  If there exists a Baby Monster that gives rise to this $T_{2}$ critical plane it must have both $v_{2} \not \equiv 0$ and 
$u_{4} \not \equiv 0$ in a neighborhood of $p_{4}$ on this submanifold.  Going one level down we have $\Delta_{3}$ is coframed by $[dv_{2} : du_{3} : dv_{3}]$ and since $u_{4} = \frac{du_{3}}{dv_{2}}$ 
and $v_{4} = \frac{dv_{3}}{dv_{2}}$ it implies that our Baby Monster must have $u_{3} \not \equiv 0$ in a neighborhood of $p_{4}$, but $v_{3}$ will vanish.  Moving another level down and looking at the
$KR$-coordinates for $\Delta_{2}$ gives $[du : du_{2} : dv_{2}]$ with $u_{3} = \frac{du}{dv_{2}}$ and $v_{3} = \frac{du_{2}}{dv_{2}}$ it implies that our Baby Monster will have $u \not \equiv 0$ and
$u_{2} \equiv 0$.  Then moving down one last level we have $\Delta_{1}$ is coframed by $[dx : du : dv]$ with $u_{2} = \frac{dx}{du}$ and $v_{2} = \frac{dv}{du}$, since $u_{2}$ is trivial it forces $x$ to be
trivial on our Baby Monster, but since $v_{2} \not \equiv 0$ it implies $v \not \equiv 0$ can't be trivial as well.  This means that both the coordinates $u \not \equiv 0$ and $v \not \equiv 0$ in a 
neighborhood of $p_{4}$ and hence we do not need to go any father down in the tower and conclude that our $T_{2}$ critical plane will come from prolonging the fiber space 
$F_{1}(p_{1}) = \pi ^{-1}_{1,0}(p_{0}) = (0,0,0, u, v)$ and that we will end up with $T_{2} = \delta ^{3}_{1}(p_{4})$ in $\Delta_{4}(p_{4})$.  In $KR$-coordinates we will have that
$\delta^{3}_{1}$ is given by $span \{ \frac{\pa}{\pa v_{2}} + u_{4} \frac{\pa}{\pa u_{3}}, \frac{\pa}{\pa u_{4}} \}$.
\

\subsubsection{Showing that $T$ will not appear in $\Delta_{4}$.}
Now, let's show that the $T$ critical plane will not appear in $\Delta_{4}(p_{4})$.  Suppose that the $T$ plane did appear in $\Delta_{4}(p_{4})$,
then in our coframing for the distribution it would be of the form $[dv_{2} : du_{4} : dv_{4}] = [a : 0 : b]$ with $a, b \in \R$ with $a \neq 0$.  Then one level down we would have in the coframing for $
\Delta_{3}$, in neighborhood of $p_{3}$, that $[dv_{2} : du_{3} : dv_{3}]$ with $u_{4} = \frac{du_{3}}{dv_{2}}$ and $v_{4} = \frac{dv_{3}}{dv_{2}}$ which implies that Baby Monster must be of the form
$[dv_{2} : du_{3} : dv_{3}] = [a : 0 : b]$ with $a, b \in \R$ and $a \neq 0$.  However we would have $u_{4} = 0$ at $p_{4}$, which creates a contradiction to what was established above about the
coordinate $u_{4} \neq 0$ at the point $p_{4}$.
\

\begin{figure}
 \label{fig:firstT2}
  \centering
  \def\svgwidth{.42 \columnwidth}
   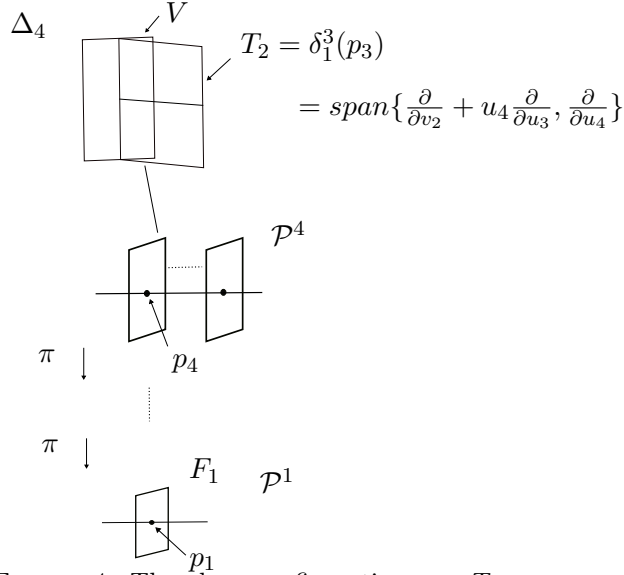
    \caption{The plane configuration over $T_{2}$.}
\end{figure}

In summary, there are two critical planes over the $T_{2}$ point: one that is the vertical plane and the other comes from a vertical plane, prolonged $3$ times, that results in a $T_{2}$
critical plane.  In Figure $4$ we see how fibers are prolonged to give us our critical plane $T_{2}$.

\begin{rem}[Continuing in the $T_{2}$ direction]
From using the above techniques in Subsection \ref{subsec:exampleT2}, it is not hard to show that if we continue moving in the $T_{2}$ directions and look at the $RVT$ code $RVLT_{2} \cdots T_{2}$
that we will end up with the exact same critical plane configuration and that the $T_{2}$ critical plane will come from the same source, the fiber space $F_{1}(p_{1})$.  This is illustrated
in Figure $6$ where the coordinates $u_{4}, \cdots, u_{k}$ cannot be zero when we evaluate at the point $p_{k}$.
\end{rem}

\subsection{The Critical Planes over the $L_{2}$ Direction.}
In this section we will show what the critical planes are over the $L_{2}$ directions for the case of $RVLL_{2}$.  Notice that this is the very first time that the $L_{2}$ direction can occur within any
$RVT$ code.  We will show that the two tangency planes will be given by
$T = \delta^{2}_{2}(p_{4})$ and $T_{2} = \delta^{3}_{1}(p_{4})$.  Initially we thought that the critical plane configuration above a $L_{2}$ point should be the same as a $T$ point or maybe a $T_{2}$ 
point, since it is contained in the intersection of those two planes.  It is worth noting that because $L_{2} = span \{ \frac{\pa}{\pa v_{2}} \}$ and when we prolong to the fourth level we have
$[dv_{2}  : du_{3} : dv_{3} ] = [1 : \frac{du_{3}}{dv_{2}} : \frac{dv_{3}}{dv_{2}} ] = [1 : 0 : 0]$ with $u_{4} = 0$ and $v_{4} = 0$ at the point $p_{4}$.  This is different than the case above for the
$T_{2}$ plane.  In that case the coordinate $u_{4} \neq 0$ and resulted in ``blocking'' the $T_{1}$ critical plane from appearing in the distribution.
\


We begin by looking at what the conditions are at each level with the $KR$-coordinates in relationship to the $RVT$ code.  Let $p_{4} \in RVLL_{2}$.
The first level: $[dx : dy : dz] = [1 : \frac{dy}{dx} :  \frac{dz}{dx}]$ to give an $R$ direction.  The second level: $[dx : du : dv] = [\frac{dx}{du} : 1 : \frac{dv}{du}]$ to give a $V$ direction.
The third level: $[du : du_{2} : dv_{2}] = [\frac{du}{dv_{2}} : \frac{du_{2}}{dv_{2}} : 1]$ to give an $L$ direction.  The $L_{2}$ direction in $\Delta_{3}(p_{3})$ is given by 
$L_{2} = span\{ \frac{\pa}{\pa v_{2}} \}$.  The fourth level is then given by $[dv_{2}: du_{3} : dv_{3}] = [1 : \frac{du_{3}}{dv_{2}} : \frac{dv_{3}}{dv_{2}}]$.  Now, $\Delta_{4}$ in
a neighborhood of $p_{4}$ will be coframed by $dv_{2}$, $du_{4}$, and $dv_{4}$.  Let's now look at what the possible critical planes are over the $L_{2}$ point.  We have
$V = span \{ \frac{\pa}{\pa u_{4}}, \frac{\pa}{\pa v_{4}} \}$.  Now, if a $T_{1}$ tangency plane is going to be present, then it will be of the form $[dv_{2} : du_{4} : dv_{4}] = [a : 0 : b]$ and
if a $T_{2}$ tangency plane is going to be present then it will be of the form $[dv_{2} : du_{4} : dv_{4}] = [a : b : 0]$.  We will work backwards to see which fiber produces these tangency critical
planes.  If we can't find such a fiber, meaning one doesn't exist, or a contradiction occurs, then no such tangency plane will exist.

\begin{figure}
  \centering
   \label{fig:prolongationfirstT2} \def\svgwidth{.65 \columnwidth}
    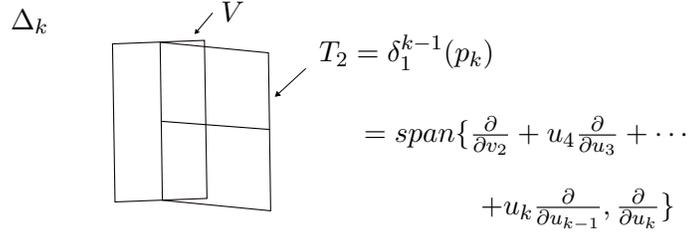
     \caption{Critical plane configuration over $p \in RVLT_{2} \cdots T_{2}$.}
\end{figure}

\subsubsection{Showing $T_{1} = \delta^{2}_{2}(p_{4})$}
We said if the $T_{1}$ direction exists then it will be given in $KR$-coordinates by $[dv_{2} : du_{4} : dv_{4}] = [a : 0 : b]$.  In local coordinates we
know, at least in a neighborhood of $p_{4}$, that the $v_{2}$ and $v_{4}$ coordinate will not be identically zero, but $u_{4}$ must be.  The tells us the following information about the local coordinate
structure of the Baby Monster.  It should be of the form: $(x,y,z, \cdots , u_{4}, v_{4}) = (0,0,0, \_, \cdots,  \_ , v_{2}, \_ , \_ ,  0 , v_{4} )$, where $v_{4}$ will be zero at $p_{4}$.  The blanks mean that
we don't have enough information to determine if these coordinates are zero or nonzero yet in a neighborhood of $p_{4}$.  Working backwards we have that $u_{4} = \frac{du_{3}}{dv_{2}}$ and
$v_{4} = \frac{dv_{3}}{dv_{2}}$.  Since $u_{4} = 0$ in a neighborhood of the point $p_{4}$ then $u_{3}$ must be as well, but since $v_{4}$ is not identically zero in a neighborhood of $p_{4}$ we must
have $v_{3}$ and $v_{2}$ are not identically zero as well.  This tells us more information about the coordinate structure, so we can fill in some of the above blanks to give
$(0,0,0, \_,  \cdots, \_ , v_{2}, 0, v_{3}, 0, v_{4})$, where both $v_{3}$ and $v_{4}$ will be zero at $p_{4}$.  Then, at the next level down we had $[du : du_{2} : dv_{2}]$ with
$u_{3} = \frac{du}{dv_{2}}$ and $v_{3} = \frac{du_{2}}{dv_{2}}$.  Since $u_{3} = 0$ in a neighborhood of the point, then $u = 0$ will be as well, but $v_{3}$ is not identically zero in a neighborhood
and means that $u_{2}$ and $v_{2}$ are not identically zero in a neighborhood of $p_{4}$ and we are able to control what values they take.  As a result, we end up with the coordinate structure
$(0,0,0, \_ , \_ ,  u_{2}, v_{2}, 0, v_{3}, 0, v_{4} )$.  This information then points to the fact that the $T_{1}$ critical
plane comes from the fiber $F_{2}(p_{2}) = (p_{1}, u_{2}, v_{2})$, where $p_{1}$ is fixed and hence $T_{1}= \delta^{2}_{2}(p_{4})$.

\subsubsection{Showing $T_{2} = \delta^{3}_{1}(p_{4})$}
The $T_{2}$ tangency plane, if it is going to appear, must be of the form $[dv_{2}, du_{4} : dv_{4}] = [a : b : 0]$.  Using the same techniques as the above, we see that in local coordinates, in a
neighborhood of $p_{4}$, it is given by $(x,y,z, \cdots, u_{4}, v_{4}) = (0,0,0, \_ ,  \cdots, \_ , v_{2}, \_ , \_ , \_ , u_{4}, 0)$.  Then one level down we must have
$[dv_{2} : du_{3} : dv_{3} ]$ and since $u_{4} = \frac{du_{3}}{dv_{2}}$ and $v_{4} = \frac{dv_{3}}{dv_{2}}$, we have the $u_{4}$ is not identically zero and so it implies that both
$u_{3}$ and $v_{2}$ can't be identically zero as well, but since $v_{4} = 0$ in a neighborhood of $p_{4}$ it implies $v_{3}$ will be as well.  This gives
$(0,0,0, \_,  \cdots , \_ , v_{2}, u_{3}, 0 , u_{4}, 0)$.  The next level down gives $[du : du_{2} : dv_{2}]$ with $u_{3} = \frac{du}{dv_{2}}$ and $v_{3} = \frac{du_{2}}{dv_{2}}$, where $u_{3}$ is not identically
equal to zero and implies that $u$ and $v_{2}$ are not as well.  Then because $v_{3} = 0$ in a neighborhood of $p_{4}$ we have the $u_{2} = 0$ in a neighborhood as well.  We now have the following 
form $(0,0,0, u, \_ , 0, v_{2}, u_{3}, 0, u_{4}, 0)$.  The next level down gives $[dx : du : dv]$ with $u_{2} = \frac{dx}{du}$, $v_{2} = \frac{dv}{du}$, where $u_{2} = 0$ in a neighborhood of $p_{4}$ 
and implies $x = 0$, but since $v_{2}$ is not identically equal zero we have $v$ is not identically zero as well.  This tells us that the coordinates are given by
$(0,0,0, u, v, 0, v_{2}, u_{3}, 0, u_{4}, 0)$, which is enough information and tells us we can stop and say that the $T_{2}$ critical plane in $\Delta_{4}(p_{4})$ will come from the prolongation of
the $F_{1}(p_{1}) = (0,0,0, u, v)$ fiber space.

\section{Proofs}
Now we are in a position to prove the two theorems of this paper.  In this section we will first present the proof of Theorem \ref{thm:spelling}, but we will only present it for two general cases.  In particular,
we will only present the cases of $R \omega L_{3}$ and $R \omega T_{2}$.  This is due to the fact that the techniques and procedure for the way the other cases are proved in an almost identical manner 
and the cases of $R \omega L_{3}$ and $R \omega T_{2}$ are illustrative examples.  Once we prove Theorem \ref{thm:spelling} we move on to Theorem \ref{thm:compspelling}.  Similar to 
the presentation of Theorem \ref{thm:spelling}, we will again only present a few illustrative examples that demonstrate the techniques used to prove the rest of the cases.  Specifically we will show
the cases $L_{2}$ and $T_{2}$.

\subsection{The Proof of Theorem \ref{thm:spelling} }
\label{sec:spelling}

In the various examples above we have already introduced the techniques and algorithm that we will use to prove this
theorem.  We will provide part of the details of the proof for the cases of
$L_{3}$ and $T_{2}$ points.
\

In addition, while all of the sources for the various critical planes were found using our algorithm of working backwards, some of the cases below like with the class $R \omega L_{2}T_{2}$ in
Subsection \ref{subsub:secLsecT}, are easier to present if we work from the level of the fiber space and then prolong them to higher levels of the Monster Tower.  This is done to help the reader follow
along with the subtleties of the $KR$-coordinate system is built, which might be somewhat difficult to derive going from the last letter of a particular $RVT$ code and then
working backwards.

\subsection{Over an $L_{3}$ point.}
\no In this section we will show that the $L_{3}$ direction has the same behavior as an $L$ direction, meaning that the plane configuration over an $L_{3}$ point will look the same as an $L$ 
direction.  While the plane configuration will look like an $L$ direction, the sources of the critical planes will be different and will depend on which letters came before the 
$L_{3}$ direction.  Let $p_{k} \in R \omega L_{3}$ and let this be the first occurrence of the $L_{3}$ 
direction in the $RVT$ code.  This implies that we can say that the code actually has to be of the form $R \omega L L_{3}$ \footnote{We should technically write $R\omega LL_{3}$ as
$R\omega' LL_{3}$, but by an abuse of notation we will just keep using $\omega$ to signify an undetermined block of $RVT$ code.  We will also continue to do this in subsequent cases.}.  In fact we can
say a little more because of the $RVT$ spelling rules
we can say that
$R \omega LL_{3} =
\begin{cases}
R \omega VT^{m}LL_{3} & \, \, \text{for} \, \, m \geq 0 \\
R \omega LLL_{3} &
\end{cases}
$
\

Now, from looking at Table \ref{tab:codes} one can see that the last $4$ cases involve the letter $T_{2}$ (these are $R \omega' LLT_{2}L_{3}$, $\cdots$, $R \omega ' V T^{m}L^{s}L_{2}T_{2}L_{3}$ for
$m \geq 0$, $s \geq 0$).  For the sake of clarity we will only present the above situation and leave these calculations to the reader, since the process used in those cases will be very similar to the 
following presentation.
\

We begin by showing that there is always at least one tangency plane over the point $L_{3}$ which is the
$\delta^{1}_{k-1}(p_{k})$ for $p_{k} \in R \omega LL_{3}$.  We start with the fiber space $F_{k-1}(p_{k-1}) = (p_{k-2}, u_{k-1}, v_{k-1})$ and use prolongation to show that it will give the
$T = \delta^{1}_{k-1}(p_{k})$ critical plane.  The first prolongation of the fiber space is given by
\begin{align*}
\mathcal{P}^{1}(F_{k-1}(p_{k-1})) & = (p_{k-2}, u_{k-1}, v_{k-1}, [dv_{k-2}: du_{k-1} : dv_{k-1}]) \\
                                                     & = (p_{k-2}, u_{k-1}, v_{k-1}, [0 : a : b]) = (p_{k-2}, u_{k-1}, v_{k-1}, [0 : 1 : \tfrac{b}{a}]) \\
                                                     & = (p_{k-2}, u_{k-1}, v_{k-1}, 0, v_{k}).
\end{align*}
Then when we evaluate at $p_{k}$ we have that $\delta ^{1}_{k-1}(p_{k})$ will be
given in the coframing for $\Delta_{k}(p_{k})$ as $[du_{k-1} : du_{k} : dv_{k} ] = [a : 0 : b]$, which gives us the $T$ critical plane.
\

Now, a second tangency plane may show up depending on the previous letters that appear in the $RVT$ code
$R \omega LL_{3}$.
\


We show that there will be at least be two tangency planes in the distribution over $L_{3}$ in $R \omega VT^{m}LL_{3}$ for $m \geq 0$.  Let $\omega$ be an $RVT$ block of length $r$.
We will show that the $T_{2}$ plane will be given by $\delta^{m+3}_{r+1}(p_{r+m+4})$ within $\Delta_{r+m+4}(p_{r+m+4})$.
Again, we will suppose that such a critical plane exists and work backwards until we arrive at the fiber space where it was born.  We begin with the coframing
$[du_{r+m+3} : du_{r + m + 4} : dv_{r + m + 4}]$ in a neighborhood of $p_{k}$ for $\Delta_{r + m + 4}$ and see that the $T_{2}$ critical plane will be of the form $[a:b:0]$ with $a \neq 0$.  Moving one 
level down we have that since $L_{3} = span \{ \frac{\pa}{\pa u_{r+m+3}} \}$, $\Delta_{r+m+3}(p_{r + m + 3})$ is coframed by
$[dv_{r+m+2} : du_{r + m + 3} : dv_{r + m + 3}]$ and our Baby Monster would be given by $[a: b: 0]$ in this coframing, telling us that $v_{r+m+2}$ and $u_{r+m+2}$ will not vanish in a neighborhood of
the point $p_{r+m+2}$ on our Baby Monster, but $v_{r+m+2}$ must vanish.  One more level down the Baby Monster in our coframing for $\Delta_{r+m+2}$ is given by
$[du_{r+1} : du_{r+m+2} : dv_{r+m+2}] = [a : 0 : b]$ since we move in an $L$ direction with $u_{r + m +3} = \frac{du_{r+1}}{dv_{r+m+2}}$ and $v_{r+m+3} = \frac{du_{r+m+2}}{dv_{r+m+2}}$.  Then after
this level we go backwards $m$ steps, where at each of those levels we moved in $T$ directions.  This gives that our Baby Monster in $\Delta_{r+2}(p_{r+2})$ will be given by
$[du_{r+1} : du_{r+2} : dv_{r+2}] = [a : 0 : b]$.  Moving one last level down we have our that our Baby Monster in $\Delta_{r+1}(p_{r+1})$ is coframed by
$[df_{r+1} : du_{r+1} dv_{r+1}] = [0 : a : b]$ since we moved in a $V$ direction, which gives $u_{r+2} = \frac{df_{r+1}}{du_{r+1}}$ and $v_{r+2} = \frac{dv_{r+1}}{du_{r+1}}$, which means we can stop at
this point and see that our $T_{2}$ critical plane comes from the fiber space $F_{r+1}(p_{r+1}) = (p_{r}, u_{r+1}, v_{r+1})$ and at level
$k = r+m+4$ is given by $\delta ^{m+3}_{r+1}(p_{r+m+4})$.
\

When we look at the $R \omega LLL_{3}$ case, we can show, using a similar argument, that for $p_{k} \in R \omega LLL_{3}$ the other critical plane comes from the prolongation of
the fiber space $F_{k-3}(p_{k-3})$.  Prolonging $3$ times gives $\delta ^{3}_{k-3}(p_{k}) = span \{ \frac{\pa}{\pa u_{k-1}}, \frac{\pa}{\pa u_{k}} \}$, which gives the $T_{2}$ critical plane in
$\Delta_{k}(p_{k})$.

\subsection{Over a $T_{2}$ point.}
\

In this section we show that there are two critical planes over $T_{2}$ points.  These are given by the vertical plane and the $T_{2}$ tangency plane.
Using $R \omega T_{2}$ as a starting point, we determine the origin of the $T_{2}$ critical plane.  This will then help us see where the plane is
situated and how it is prolonged to higher levels of the tower.  Based on our knowledge of the spelling rules, and the fact that we
are assuming that no other $T_{2}$ occurs in $R \omega T_{2}$, the $RVT$ codes look like either $R \omega L T_{2}$, $R \omega L_{2} T_{2}$, or
$R \omega L_{3} T_{2}$, which is based upon the fact that we are assuming only one occurrence of each of the $L_{j}$'s for $j = 1,2,3$ occurring in the code as well.  Start with the case of $R \omega L
T_{2}$.  We have
\[ R \omega L T_{2} =
\begin{cases}
R \omega V T^{m}LT_{2} & \, \, \text{for} \, \, m \geq 0 \\
R \omega L L T_{2} &
\end{cases}
\]

\

\no For the classes of the form $R \omega L_{2}T_{2}$ there will be two possible forms, where
\

\no \[
R \omega L_{2}T_{2} =
\begin{cases}
R \omega VT^{m}L L_{2} T_{2} & \, \, \text{for} \, \, m \geq 0 \\
R \omega V T^{m} L^{s} L_{2} T_{2} & \, \, \text{for} \, \, m \geq 0, s \geq 2
\end{cases}
\]

\

\no Classes of the form $R \omega L_{3}T_{2}$ will be of the following form
\

\no \[
R \omega L_{3}T_{2} =
\begin{cases}
R \omega V T^{m}LL_{3}T_{2} & \, \, \text{for} \, \, m \geq 0 \\
R \omega VT^{m}L^{s}L_{3}T_{2} &  \, \, \text{for} \, \, m \geq 0, s \geq 2
\end{cases}
\]

\subsubsection{The class $R \omega T_{2} = R \omega L T_{2}$ }
\label{subsub:secLsecT}
\

We start with $R \omega V T^{m}LT_{2}$ with $m \geq 0$.  Before we begin we want to point out that we can actually look just at the $RVT$ code $RVT^{m}LT_{2}$ instead of $R \omega VT^{m}LT_{2}$ because of which fiber space we need to prolong in order
to obtain the $T_{2}$ critical plane.  In short, we do not need to go any farther back than the level where the last $V$ appears.  This will also ease the amount of notation needed.
\

We want to show for the code $RVT^{m}LT_{2}$ for $m \geq 0$, that $T_{2}$ comes from the critical plane given by
$T_{2} = \delta^{m+3}_{1}(p_{m+4})$.  We begin by looking in a neighborhood of the point $p_{m+4}$ and we have that $T_{2}$, if it exists, in $\Delta_{m+4}$ will be of the form
$[dv_{m+2} : du_{m+4} : dv_{m+4}] = [a : b : 0]$ with $a \neq 0$.  Two levels down we have $\delta ^{m+2}_{1}(p_{m+3})$ will be given in $\Delta_{m+3}(p_{m+3})$ by
$[du : du_{m+2} : dv_{m+2}] = [a : 0 : b] = [\frac{a}{b} : 0 : 1]$ to give a $L$ direction.
\

The fiber $F_{1}(p_{1})$ gives the vertical plane $[dx: du :dv] = [0 : a : b]$ within $\Delta _{1}(p_{1})$.  After moving in a $V$ direction the critical plane $\delta ^{1}_{1}(p_{s})$ will be
given in the coframing for $\Delta_{2}(p_{2})$ by $[du: du_{2} : dv_{2}] = [a : 0 : b] = [1 : 0 : \frac{a}{b}]$ to give a $T$ direction.  If one were to
continue to move in strictly $T$ directions $(m-1)$ more times that this pattern will persist and we will have $\delta^{m+1}_{1}(p_{m+2})$ will be given by $[du : du_{m+1} :
dv_{m+1}] = [a : 0 : b] = [\frac{a}{b} : 0 : 1]$, moving in an $L$ direction.  This shows that the $T_{2}$ critical plane will indeed come from the $F_{1}(p_{1})$ fiber.
\

The next step is to show that over $T_{2}$ the only spelling rules are $R$, $V$, and $T_{2}$.  Suppose that $\ell \subseteq \Delta_{m+3}(p_{m+3})$ and is strictly in the $T_{2}$ plane.  This
means for a direction $\ell$ to lie in this plane it will be of the form $\ell = a \frac{\pa}{\pa v_{m+2}} + b \frac{\pa}{\pa u_{m+3}}$ with $a,b \in \R$ with $a \neq 0$.  The $KR$-coordinates for the $(m+4)$-th 
level of the tower will be $[dv_{m+2} : du_{m+3} : dv_{m+3} ] = [1 : \frac{du_{m+3}}{dv_{m+2}}  :  \frac{dv_{m+3}}{dv_{m+2}}]$.  We will have $u_{m+4} = \frac{du_{m+3}}{dv_{m
+2}}$ and $v_{m+4} = \frac{dv_{m+3}}{dv_{m+2}}$ where $u_{m+4}$ cannot be equal to zero at the point $p_{m+4}$ or else it will not be a $T_{2}$ point.
\

In summary, the above argument shows that the $T$ critical plane will not appear in the distribution over points $p_{k} \in RVT^{m}LT_{2}$, since it would force $u_{m+4} = 0$ at the point $p_{m+4}$, and 
hence in the distribution over points $p_{k} \in R \omega VT^{m}LT_{2}$.
We would like to point out that the above argument can be repeated in a very similar fashion to show that the $T_{1}$ critical planes will not appear over any of these points.  In order to avoid being
repetitive we will omit presenting this justification for the absence of the $T$ critical plane again and again for each of the following cases below.

\begin{rem}[Continuing to Move in $T_{2}$ Directions]
From the code $R VT^{m}LT_{2}$ it can be shown that if we continue to move in direction within the plane $T_{2}$ then the $RVT$ code will have the form $RVT^{m}LT^{s}_{2}$ and 
the $T_{2}$ will have the form $\delta^{m+s+1}_{1}(p_{m+s+2}) = span \{ \frac{\pa}{\pa v_{m+2}} + u_{m+5} \frac{\pa}{\pa u_{m+4}} + \cdots + u_{m+s+2} \frac{\pa}{\pa u_{m+s+1}} \frac{\pa}{\pa u_{m+s
+2}} \}$.
\end{rem}

Now we present the case of $R \omega ' LL T_{2}$.  We assume the total length of this code is $k$ and the coframing for the the distribution $\Delta_{k}$ in a neighborhood of points in $R \omega ' LL 
T_{2}$ is given in $KR$-coordinates by $[dv_{k-2} : du_{k} : dv_{k}]$ and if the $T_{2}$ plane exists it is given by $[a : b : 0]$ in this coframing.  We aim to show that in $\Delta_{k}$ that 
$T_{2} = \delta^{3}_{k-3}$.  Going one level down, this plane projects to $[a : b : 0]$ in the coframing $[dv_{k-2}: du_{k-1}:dv_{k-1}]$ for $\Delta_{k-1}$.  Again, projecting this plane another level
down to $\Delta_{k-2}$ we have $[dv_{k-3} : du_{k-2} : dv_{k-2}] = [a : 0 : b]$.  Doing this one last time, the plane is given in the coframing for $\Delta_{k-3}$ by $[df_{k-4} : du_{k-3} : dv_{k-3}] = 
[0 : a : b]$, which indicates that the $T_{2}$ plane in $\Delta_{k}$ does in fact come from the fiber space $F_{k-3}$ and is equal to $\delta^{3}_{k-3}$.  We can also show that the $T$ plane does
not exist in $\Delta_{k}$ by using a similar argument as the above case for $R \omega V T^{m}LT_{2}$ because the coordinate $u_{k} \neq 0$ at any point in $R \omega ' LL T_{2}$.

\subsubsection{The class $R \omega T_{2} = R \omega L_{2} T_{2}$ }
\label{subsub:secL2secT}
\

In this case we will start with $R \omega V T^{m} L L_{2} T_{2}$.  We will again, for ease of notation, just look at the code $RVT^{m} L L_{2} T_{2}$, where one can see how it would generalize
for the case of $R \omega V T^{m} L L_{2} T_{2}$.   We show that the $T_{2}$ critical plane that appears over the points in $RVT^{m} L L_{2} T_{2}$ will come from the Baby Monster
$\delta ^{m+4}_{1}$.  We begin with the fiber space $F_{1}(p_{1})$ and see that the vertical plane $\delta ^{0}_{1}(p_{1})$ will be given in $\Delta_{1}(p_{1})$ as
$[dx : du : dv] = [0 : a : b]$ for $a, b \in \R$ with $a \neq 0$.  Then $\delta ^{1}_{1}$ will be given in a neighborhood of the point $p_{2}$ in $\Delta_{2}$ by $[du : du_{2} : dv_{2}] =
[a : 0 : b] = [1 : 0 : \frac{b}{a}]$ to give a $T$ direction.  Repeating this $m-1$ times more gives $\delta ^{m}_{1}$ in a neighborhood of $p_{m+1}$ by
$[du : du_{m+1} : dv_{m+1}] = [a : 0 : b] = [1 : 0 : \frac{b}{a}]$ to give the last $T$ direction.  For the $L$ direction we have that $\delta ^{m+1}_{1}$ will be $[du : du_{m+2} : dv_{m+2}] =
[a : 0 : b] = [\frac{a}{b} : 0 : 1]$.  For the $L_{2}$ direction we have that $\delta ^{m+2}_{1}$ will be $[dv_{m+2} : du_{m+3} : dv_{m+3}] = [a : b : 0] = [1 : \frac{b}{a} : 0]$.  Then for $\delta ^{m+3}_{1}$ we 
have $[dv_{m+2} : du_{m+4} : dv_{m+4}] = [a : b : 0] = [1 : \frac{b}{a} : 0]$ to give the $T_{2}$ direction.  Lastly, $\delta ^{m+4}_{1}$ will be given in $\Delta_{m+5}$ within a 
neighborhood of $p_{m + 4}$ by $[dv_{m+2} : du_{m+5} : dv_{m+5}] = [a : b : 0] $, which gives the desired result.
\

The next case is when $R \omega L_{2}T_{2} = R \omega V T^{m}L^{s}L_{2} T_{2}$ for $m \geq 0$ and $s \geq 2$.  We will look at the $RVT$ code $R V T^{m}L^{s}L_{2}T_{2}$ and show that the
$T_{2}$ critical plane in $\Delta_{m+s+4}(p_{m+s+4})$ is given by the Baby Monster $\delta ^{4}_{m + s} (p_{m+s+4})$.  Start with local coordinates for $\Delta_{m+s}$ in a neighborhood of
$p_{m+s} \in RVT^{m}L^{s-2}$ for $p_{m+s+4} \in RVT^{m} L^{s}L_{2}T_{2}$, which is given by $[dv_{m+s-1} : du_{m+s} : dv_{m+2}]$ and when we prolong the fiber space $F_{m+s}(p_{m+s})$ we will 
get in a neighborhood of $p_{m+s+1}$ that $\delta ^{1}_{m+s}$ is $[dv_{m+s} : du_{m+s+1} : dv_{m+s+1}] = [a : 0 : b]$.  After going in an $L$ direction we have locally that $\delta ^{2}_{m+s}$ is given by
$[dv_{m+s+1} : du_{m+s+2} : dv_{m+s+2}] = [a : b : 0]$.  Then after going in the $L_{2} = span \{ \frac{\pa}{\pa v_{m+s+1}} \}$ direction, $\delta ^{3}_{m+2}$ is given by
$[ dv_{m+s+1} : du_{m+s+3} : dv_{m+s+3}] = [ a : b : 0]$ and after moving in a $T_{2}$ direction, where we divide by $dv_{m+s+1}$, we get $\delta ^{4}_{m+s}$ in a neighborhood of $p_{m+s+4}$ is
given in the $KR$-coordinates for $\Delta_{m+s+4}$ by $[dv_{m+s} : du_{m+s+4}  : dv_{m+s+4} ] = [a : b : 0]$, which gives us the desired $T_{2}$ critical plane.

\subsubsection{The class $R \omega L_{3}T_{2}$}
\

The argument for this case is similar to the previous two cases presented in Sections \ref{subsub:secLsecT} and \ref{subsub:secL2secT} above.  As a result, we will omit this case.

\subsection{The Proof of Theorem \ref{thm:compspelling} }

In this section we will show using induction how the spelling rules hold in the general case.  We understand how the spelling rules work for the cases of $R, V, T,$ and $L$.  We want to
completely determine the spelling rules for the letters $T_{2}$, $L_{2}$, and $L_{3}$.  We have already shown what the spelling rules are for the first occurrence of each of these letters
in a given $RVT$ code in Table $2$.  We need to take an arbitrary $RVT$ code that ends in either a $T_{2}$, $L_{2}$, or $L_{3}$ and show that these spelling rules will persist.  
We do this by induction on the number of $T_{2}$'s, $L_{2}$'s, and $L_{3}$'s in a given $RVT$ code.  We pause for a moment to explain this method of induction on the total number of 
$T_{2}$'s, $L_{2}$'s, and $L_{3}$'s through a few examples.  If $n=4$, then $RVLL_{2}T_{2}L_{3}L_{2}$, $RRVL^{4}_{2}$, $RVLT_{2}T_{2}L_{2}L_{2}$, and
$RVT^{7}LL_{3}L_{2}L_{3}L_{2}$ are examples.  However, $RVLL_{2}L_{3}T_{2}T_{2}L_{3}$ and $RVTLL_{2}T_{2}R$ are not examples, since there is a total of $5$ $T_{2}$'s, $L_{2}$'s, and $L_{3}$'s
in the first example and the sec on example only has a total of $2$ of these letters.  
\

The base case is $n=1$, proved above in Section \ref{sec:spelling}, and is detailed in Table $2$.  
The inductive hypothesis is that the spelling rules hold true for any $RVT$ code that has a total number of $n$ $T_{2}$'s, $L_{2}$'s, or $L_{3}$'s in it.  Now suppose we have an $RVT$ code $\omega$ 
that contains a total of $n+1$ of these letters in it and suppose that the last letter of $\omega$ is one of these three possible letters.  We show that the spelling rules will hold for the cases of $\omega = 
\omega' L_{2}$ and $\omega = \omega' T_{2}$, where the terminal letter in the code is $L_{2}$ and $T_{2}$ respectively, from which one can see how the case of $\omega = \omega' L_{3}$ will follow.

\subsection{The case of the terminal letter $L_{2}$}
\

Our code is $\omega' L_{2}$ and we want to show that the letters $R$, $V$, $T = T_{i}$ for $i = 1,2$ and $L_{j}$ for $j = 1,2,3$ exist over this $L_{2}$.  This means we need to look at the possible
planes that appear in $\Delta_{k} (q)$, where $k$ equals the length of the code $\omega$ and $q$ is a point in a neighborhood of $p \in \omega$.  By the inductive hypothesis we have the following
refinement of our code
\[ 
\omega' L_{2} =
\begin{cases}
\tilde{\omega} L_{1} L_{2} \\
\tilde{\omega} L_{2} L_{2} \\
\tilde{\omega} L_{3} L_{2}
\end{cases}
\]

\no A coframing for our distribution $\Delta_{k}(q)$ is given by $[df_{k-1} : du_{k} : dv_{k}]$.  Also, we will show that the spelling rules hold for the case of $\tilde{\omega} L L_{2}$,
since the cases of $\tilde{\omega} L_{j} L_{2}$ for $j = 2,3$ are almost identical.

\subsubsection{Case of $\tilde{\omega}L L_{2}$}
\

We can refine our coframing and say that our uniformizing coordinate is $df_{k-1} = dv_{k-2}$.  This is because in $\Delta_{k-2}$ we move in an $L$ direction to give
$[df_{k-3} : du_{k-2} : dv_{k-2}] = [ \frac{df_{k-3}}{dv_{k-2}} : \frac{du_{k-2}}{dv_{k-2}}  : 1]$ and then in a $L_{2}$ direction $[dv_{k-2} : du_{k-1} : dv_{k-1}] =
[1 : \frac{du_{k-1}}{dv_{k-2}} : \frac{dv_{k-1}}{dv_{k-2}}]$ in $\Delta_{k-1}$.
\

Now, if the $T$ critical plane exists in $\Delta_{k}(q)$
then it must be of the form $[dv_{k-2}: du_{k}: dv_{k}] = [a: 0 : b]$.  When we project one level down we have $[dv_{k-2}: du_{n-1}: dv_{n-1}] =
[dv_{k-2} : du_{k-1} : dv_{k-1}] = [a : 0 : b]$ in $\Delta_{k-1}$ which is the $T$ plane.  This implies that the $T$ critical plane will exist in $\Delta_{k}(q)$ and comes from the
prolongation of the $T$ plane one level below.
\

If the $T_{2}$ critical plane is to exist in $\Delta_{k}(q)$ then it is of the form $[dv_{k-2}: du_{n}: dv_{n}] = [a : b : 0]$.  Projecting one level down we have
$[dv_{k-2} : du_{k-1} : dv_{k-1}] = [a : b : 0]$, which is the $T_{2}$ critical plane, which must exist in $\Delta_{k-1}(q)$ since we are over an $L$ point.
\

Since both tangency critical planes exist above $\omega L L_{2}$ we end up with the distinguished directions $L_{j}$ for $j = 1,2,3$ existing as well in $\Delta_{k}(q)$ and hence the spelling
rules persist.

\subsection{The case of the terminal letter $T_{2}$}
\

We want to show that from the code $\omega' T_{2}$ that the letters $R, V, T_{2}$, and $L_{3}$ are the only possibilities.  We again need to look at the possible planes that appear
in $\Delta_{k}(q)$.  By the inductive hypothesis we have the following refinement of our code
$\omega ' T_{2} =
\begin{cases}
\tilde{\omega} L T_{2} \\
\tilde{\omega} L_{2} T_{2} \\
\tilde{\omega} L_{3} T_{2} \\
\tilde{\omega} T_{2} T_{2}
\end{cases}
$
We will show that the spelling rules hold for $\tilde{\omega} L T_{2}$ and $\tilde{\omega} T_{2}T_{2}$, since the cases of $\tilde{\omega} L_{j} T_{2}$ for $j = 2,3$ will be similar.

\subsubsection{Case of $\tilde{\omega}LT_{2}$}
\

We look at the code $\tilde{\omega} L T_{2}$ of length $k$, where $\Delta_{k}$ is coframed by $[df_{k-1} : du_{k} : dv_{k}] = [dv_{k-2} : du_{k} : dv_{k}]$.  This is because when we go two levels down
we have that $\Delta_{k-2}$ is coframed by $[df_{k-3} : du_{k-2} : dv_{k-2}] = [\frac{df_{k-3}}{dv_{k-2}} : \frac{du_{k-2}}{dv_{k-2}} : 1]$ to give an $L$ direction.  Notice that when we move in a $T_{2}$ 
direction in $\Delta_{k-1}$ we have $[dv_{k-2} : du_{k-1} : dv_{k-1}] = [a : b : 0] = [1 : \frac{b}{a} : 0]$.  This results in $u_{k} \neq 0$ and $v_{k}$ is identically zero for points in $\tilde{\omega} T_{2}$, 
where $u_{k} = \frac{du_{k-1}}{dv_{k-2}}$ and $v_{k} = \frac{dv_{k-1}}{dv_{k-2}}$.
\

For $T_{2}$ to be present in $\Delta_{k}$ we would have $[dv_{k-2} : du_{k} : dv_{k}] = [a : b : 0]$ and projecting
one level down gives $[dv_{k-2} : du_{k-1} : dv_{k-1}] = [a : b : 0]$, the $T_{2}$ critical plane in
$\Delta_{k-1}$ and again we get that $T_{2}$ will exist in $\Delta_{k}$.
\

We now show that no $T$ critical plane can exist in $\Delta_{k}$.  Suppose that such a plane did exit, then it would be given by
$[dv_{k-2} : du_{k} : dv_{k}] = [a : 0 : b]$.  Projecting one level down gives $[dv_{k-2} : du_{k-1} : dv_{k-1}] = [a : 0 : b]$.  When we move in a $T_{2}$ direction
we have $[1 : 0 : \frac{b}{a}]$ and gives $u_{k} = 0$ for points within $\omega ' T_{2}$, but this creates a contradiction to $u_{k} \neq 0$.  Hence $T$ is not present in $\Delta_{k}$.
\

Then since $L_{3}$ comes from the intersection of the vertical plane and the $T_{2}$ plane, we have shown the spelling rules over the last $T_{2}$.

\subsubsection{Case of $\tilde{\omega}T_{2}T_{2}$}
\

From the code $\omega' T_{2} T_{2}$ we have the coframing $[df_{k-1} : du_{k} : dv_{k}] = [df_{k-3} : du_{k} : dv_{k}]$ for $\Delta_{k}$, since in $\Delta_{k-2}$ has
$[df_{k-3} : u_{k-2}: dv_{k-2}] = [1 : \frac{du_{k-2}}{df_{k-3}} : \frac{dv_{k-2}}{df_{k-3}}] = [1 : u_{k-1} : v_{k-1}]$ with $u_{k-1} \neq 0$ and $v_{k}$ vanishing at points in $\omega' T_{2}$ and similarly
in $\Delta_{k-1}$ we have $[df_{k-2} : du_{k-1} : dv_{k-1}] = [df_{k-3} : du_{k-1} : dv_{k-1}] = [ 1 : \frac{du_{k-1}}{df_{k-2}} : \frac{dv_{k-1}}{df_{k-2}} ] = [1 : u_{k} : v_{k}]$, with $u_{k} \neq 0$ and
$v_{k} = 0$ at points in $\omega' T_{2} T_{2}$.
\

For $T_{2}$ to be present in $\Delta_{k}$ we would have $[df_{k-3} : du_{k} : dv_{k}] = [a : b : 0]$ for $a, b \in \R$ for $a \neq 0$.  Projecting one level down gives $[df_{k-3} : du_{k-1}: dv_{k-1}] =
[a: b : 0]$, which is again the $T_{2}$ critical plane.  Therefore the $T_{2}$ critical plane will exist in $\Delta_{k}$.
\

Lastly, we show that $T_{1}$ cannot be present within the distribution $\Delta_{k}$.  Suppose that it was, then it would be given by $[df_{k-3} : du_{k} : dv_{k}] =
[a : 0 : b]$ and we project one level down to get $[df_{k-3} : du_{k-1} : dv_{k-1}] = [a : 0 : b] = [1 : 0 : \frac{b}{a}] = [1 : u_{k} : v_{k}]$, but would imply that $u_{k} = 0$, which creates a contradiction.

\section{Conclusion}

In this work we have completed the spelling rules for the $RVT$ code in the $\R ^{3}$- Monster Tower and a method which determines the level at which the critical planes are born from within the 
rank $3$ distribution.
\

It is also worth pointing out that this algorithm can be applied to a more general problem of understanding the incidence relations that occur within the $\R^{n}$-Monster Tower.  In this case $n \geq 4$, as 
was pointed out in \cite{pelletier2}, there will be even more critical planes existing, which makes the $RVT$ class spelling rules even more complicated.
The algorithm presented in this paper serves as a tool to help understand what incidence relations occur in this more general case as well.
\

One of the main motivations for this work was not only to understand the geometry of the Monster Tower, but also to understand the relationship between the $RVT$
coding system and Mormul's $EKR$ coding system.  The authors first studied this relationship in \cite{castro2} in order to answer some basic questions concerning connections
between the $RVT$ coding system and the $EKR$ system that Mormul uses to stratify Goursat Multi-Flags \cite{mormul1}.  Recent work done by Pelletier and Slayman in
\cite{pelletier2} displays not only the connection between the two coding systems up to level $4$ of the Monster Tower, but also establishes how the $EKR$ coding system
relates to the constraints that are placed on the various configurations of the articulated arm system.  Pelletier and Slayman present a beautiful result, Proposition $5.1$, that displays
the relationship between the $EKR$ codes and these restrictions placed on the articulated arms.  In addition, they also raise an important question about how certain $EKR$ classes give rise 
to more than one $RVT$ class.   An example of this can be seen with the $EKR$ class $121$ splitting into the $RVT$ classes $RVR$ and $RVT$.  Further examples can be found in \cite{castro2} and
\cite{howard}.  Pelletier and Slayman point out that it would be interesting to see if the articulated arm model could explain this splitting between the $RVT$ and $EKR$ systems.
Taking these two things into consideration, it would be interesting to use these connections to further explore the relationship between the geometry of the Monster Tower and the articulated arm system.


%



\medskip
\bibliographystyle{alpha}  
\bibliography{orbitspart2ver-4}


\end{document}